\documentclass[11pt]{article}
\usepackage{amsmath,amstext,amscd, amssymb, epsfig, psfrag}

\newtheorem{thm}{Theorem}[section]
\newtheorem{lemma}[thm]{Lemma}
\newtheorem{cor}[thm]{Corollary}
\newtheorem{prop}[thm]{Proposition}
\newtheorem{rem}[thm]{Remark}

\newtheorem{defn}[thm]{Definition}

\newtheorem{remark}[thm]{Remark}

\newtheorem{Open Questions}[thm]{Open Questions}
\newtheorem{Open Question}[thm]{Open Question}

\def\cal{\mathcal}

\def\Bbb{\mathbb}
\def\bar{\overline}
\def\sstar{\star}

\def\Z{\Bbb{Z}}
\def\N{\Bbb{N}}

\def\ni{\noindent}

\def\Cay{\hbox{\it Cay}}

\def\EDiam{\hbox{\rm EDiam}}
\def\IDiam{\hbox{\rm IDiam}}

\def\F+L{\hbox{$\textup{F}\!_+\textup{L}$}}

\def\ssm{\smallsetminus}

\def\ms{\medskip}

\def\onto{{\kern3pt\to\kern-8pt\to\kern3pt}}

\def\<{\langle}
\def\>{\rangle}
\def\|{{\ |\ }}

\def\G{\Gamma}

\def\s{\sigma}

\def\OO{\cal O}
\def\TT{\cal T}

 \def\AA{\cal A}
\def\BB{\cal B}
 \def\CC{\cal C}
 \def\RR{\cal R}
 \def\SS{\cal S}

 \def\PP{\cal P}
 \def\UU{\cal U}
\def\QQ{\cal Q}

\newcommand{\set}[1]{\left\{#1\right\}}
\newcommand{\restricted}[1]{\left|_{#1} \right.}
\newcommand{\abs}[1]{\left|#1\right|}
\renewcommand{\ni}{\noindent}
\renewcommand{\ss}{\smallskip}
\renewcommand{\ms}{\medskip}
\newcommand{\bs}{\bigskip}
\def\*{^{\sstar}}

\newcommand{\Proof}{\ni \emph{Proof. }}
\renewcommand{\t}{\tilde}
\newcommand{\B}[2]{\left(\!\!\!\begin{array}{c}{#1} \\ {#2}\end{array}\!\!\!\right)}
\newcommand{\TB}[2]{\mbox{\tiny{$\left(\!\!\!\begin{array}{c}{#1} \\
{#2}\end{array}\!\!\!\right)$}}}
\newcommand{\SB}[2]{\mbox{\footnotesize{$\left(\!\!\!\begin{array}{c}{#1} \\
{#2}\end{array}\!\!\!\right)$}}}

\def\qed{{\ifhmode\unskip\nobreak\hfil\penalty50 \hskip1em \else\nobreak\fi
   \mbox{}\nobreak\hfil \rule{1ex}{1ex}
   \parfillskip=0pt \finalhyphendemerits=0 \par}}

\begin{document}

\title{\textbf{Extrinsic versus intrinsic diameter \\ for Riemannian 
filling--discs \\ and van Kampen diagrams}}

\author{M.\ R.\ Bridson
 and T.\ R.\ Riley\footnote{The first author was supported in part by fellowships from the EPSRC of Great Britain and the second author in part by NSF grant  0404767.}}


\date{}


\maketitle

\begin{abstract} 
\noindent The diameter of a disc filling
a loop in the universal covering 
of a Riemannian manifold $M$ may be measured {\em{extrinsically}} using the
distance function on the ambient space or \emph{intrinsically} using the
induced length metric on the disc. Correspondingly, the diameter of a 
 van~Kampen diagram $\Delta$ filling a word that
represents the identity in a finitely presented group $\Gamma$ can either be measured
\emph{intrinsically} in the 1--skeleton of $\Delta$  or
\emph{extrinsically} in the Cayley graph of $\Gamma$. We construct
the first examples of closed manifolds $M$ and finitely
presented groups $\G=\pi_1M$ for which this choice --- intrinsic
versus extrinsic --- gives rise to qualitatively different min--diameter
filling functions.\ms
\\
\footnotesize{\ni \textbf{2000 Mathematics Subject
Classification: 53C23, 20F65, 20F10}  \\ \ni \emph{Key words and phrases:} diameter,
isodiametric inequality, van~Kampen diagram}
\end{abstract}

\section{Introduction and preliminaries} \label{intro}

A continuous map $F : \mathbb{D}^2 \to \widetilde{M}$ from the unit disc $\mathbb{D}^2 \subset \mathbb R^2$ to the universal cover of a Riemannian manifold $M$  is called a \emph{disc--filling} for a loop $c$ in $\widetilde{M}$ when the restriction of $F$ to $\partial \mathbb{D}^2$ is a monotone reparameterisation of $c$. According to context, one might measure the diameter of such a
disc--filling \emph{intrinsically} using the pseudo--distance
$D(a,b)$ on $\mathbb{D}^2$ defined as the infimum of the lengths of
 curves $F \circ p$ such that $p:[0,1] \to \mathbb{D}^2$ is a path
from $a$ to $b$,
 or \emph{extrinsically} using the distance function $d$ on
 the  ambient manifold:
\begin{eqnarray*}
\IDiam(F) & := & \sup  \left\{ \,  D(a,b)  \mid  a,b \in \mathbb{D}^2
\,  \right\}  \  \text{ and} \\
\EDiam(F) & := & \sup \left\{ \,  d(F(a), F(b))  \mid a,b \in
\mathbb{D}^2   \,    \right\}.
\end{eqnarray*}
   
   These two options give different functionals on the space of 
   rectifiable loops in 
$\widetilde{M}$:  the \emph{intrinsic} and \emph{extrinsic} diameter functionals $\IDiam_{\widetilde{M}}, \EDiam_{\widetilde{M}}:[0,\infty) \to[0,\infty)$, defined by 
\begin{eqnarray*}
\IDiam_{\widetilde{M}}(l) & := & \sup_{c} \,  \inf_{F}  \,  \IDiam(F) \  \text{ and} \\ 
\EDiam_{\widetilde{M}}(l) & := &  \sup_{c} \,  \inf_{F} \,  \EDiam(F) 
\end{eqnarray*}
where, in each case, the supremum ranges over loops $c$ in $\widetilde{M}$ of length at most $l$ and the infimum over disc--fillings of $c$.

 There is an obvious inequality $\IDiam(F) \geq \EDiam(F)$, which passes to the functionals:  $\IDiam_{\widetilde{M}}(l) \geq \EDiam_{\widetilde{M}}(l)$ for all $l$.  One anticipates
that for certain compact manifolds  $M$, families of minimal intrinsic diameter filling--discs might \emph{fold--back} on themselves so as to have  smaller extrinsic than intrinsic diameter, and the two functionals might then differ asymptotically.  
However, it has proved hard to animate this
intuition with examples. In this paper we overcome this difficulty.

\begin{thm} \label{main thm1}
For every $\alpha>0$, there is a closed connected Riemannian manifold $M$
and some $\beta >0$ such that $\EDiam_{\widetilde{M}}(\ell) \preceq \ell^{\beta}$ and
$\IDiam_{\widetilde{M}}(\ell) \succeq \ell^{\alpha+\beta}$.
\end{thm}

[For functions
$f,g : [0,\infty) \to [0,\infty)$  we write $f \preceq g$ when there exists
$C>0$ such that $f(l) \leq C g(Cl+C) + C l +C$ for all $l \in [0,\infty)$,  and we say $f \simeq g$ when $f \preceq g$ and $g \preceq f$.  We extend these relations to accommodate functions with domain $\N$ by declaring such functions to be constant on the half--open intervals $[n,n+1)$.] 

\ms

Our main arguments will not be cast in the language of Riemannian manifolds
but in terms of combinatorial filling problems for finitely
presented groups.  There is a close correspondence, given full voice by M.~Gromov \cite{Gromov, Gromov5}, between
the geometry of filling--discs in compact Riemannian manifolds
$M$ and the complexity of the word problem in $\pi_1M$. Most famously,
the 2--dimensional, genus--0 isoperimetric function of $M$ has the same
asymptotic behaviour as (i.e.\ is $\simeq$ equivalent to) the Dehn function
of any finite presentation 
of $\pi_1M$. (See
\cite{Bridson6} for details.) But the correspondence is wider: most natural ways of bounding the geometry of filling--discs in a manifold
correspond to measures of the complexity of the word
problem in $\pi_1M$. The following theorem, proved in Section~\ref{Riemannian-algebraic}, says that it applies to $\IDiam_{\widetilde{M}}$ and $\EDiam_{\widetilde{M}}$, and their group theoretic analogues $\IDiam_{\PP}$ and $\EDiam_{\PP}$ (as defined below).

\begin{thm} \label{Riem-Comb}
If $\PP$ is a finite presentation of the fundamental group of a closed,
 connected,  Riemannian manifold $M$ then $$\IDiam_{\PP} \simeq \IDiam_{\widetilde{M}} \ \  \textit{ and } \ \  \EDiam_{\PP}
 \simeq \EDiam_{\widetilde{M}}.$$
\end{thm}

Let $\PP=\langle \AA \mid \RR \rangle$ be a finite presentation of a group $\Gamma$. The length of a 
\emph{word} $w$ in the free
monoid $(\AA^{\pm 1})^{\ast}$   is denoted $\ell(w)$.    Let $d$ denote the word metric on $\Gamma$ associated to $\mathcal{A}$.   
If $w=1$ in $\G$ (in which case $w$ is called \emph{null--homotopic}) 
then there exist equalities in the free group $F(\AA)$ of the form
$$
w \ = \ \prod_{i=1}^N {u_i}^{-1}r_iu_i,
$$
where $r_i \in \RR^{\pm 1}$ and $u_i \in (\AA^{\pm 1})^{\ast}$.  Ranging over all such products $w_0$ freely equal to $w$, define $\IDiam_{\PP}(w)$ to be the minimum  of $\max_{1 \leq i \leq N} \ell(u_i)$, and define $\EDiam_{\PP}(w)$ to be the minimum of $$\max \set{ \,  d(1,\bar{p}) \,  \mid \, p \textup{ a prefix of } w_0  \, }$$ where $\bar{p}$ denotes the element of $\Gamma$ represented by $p$. 
Define $\IDiam_{\PP},\, \EDiam_{\PP}: \N \to \N$ by letting  $\IDiam_{\PP}(n)$ and $\EDiam_{\PP}(n)$ be the maxima of $\IDiam(w)$ and $\EDiam(w)$, respectively, over all null--homotopic words $w$ with  $\ell(w) \leq n$. 

As we will explain in Section~\ref{Riemannian-algebraic}, one can reinterpret $\IDiam_{\PP}(w)$ and $\EDiam_{\PP}(w)$ using ``combinatorial filling--discs'' (\emph{van~Kampen diagrams}) for loops in the Cayley 2--complex of $\PP$.  Then relating the Cayley 2--complex of $\PP$ to  $\widetilde{M}$ will lead to a proof of Theorem~\ref{Riem-Comb}.   

In Section~\ref{the groups} we describe an infinite family of finite presentations whose diameter functionals exhibit the range of behaviour described in the following theorem.

\begin{thm} \label{main thm}
For every $\alpha>0$, there is a group with  finite presentation $\PP$
and some $\beta >0$ such that $\EDiam_{\PP}(n) \preceq n^{\beta}$ and
$\IDiam_{\PP}(n) \succeq n^{\alpha+\beta}$.
\end{thm}

Theorem~\ref{main thm1} is an immediate consequence of 
Theorems~\ref{Riem-Comb} and \ref{main thm} because, as is well known, every finitely presentable 
group is the fundamental group of a closed connected Riemannian manifold.

\ms

Our results prompt the following questions.   What is the optimal upper bound for $\IDiam_{\PP}$ in terms of $\EDiam_{\PP}$ for general finite presentations $\PP$?  How might one construct a finite
presentation $\PP$ for which $\EDiam_{\PP}(n)$ is bounded above by a polynomial function and $\IDiam_{\PP}(n)$ is bounded below by an exponential function?

\subsection{A sketch of the proof of Theorem~\ref{main thm}}

We prove Theorem~\ref{main thm} by constructing a novel family of
groups  $\Psi_{k,m}$. They are obtained by amalgamating 
auxiliary
groups $\Gamma_m$ and $\Phi_k$ along an infinite cyclic  subgroup generated by a letter $t$.  

For each positive integer $n$, one finds an edge--loop in the
Cayley 2--complex of  $\Gamma_m$ that has length roughly
 $n$ and intrinsic diameter roughly $n^m$.  And the reason for this large intrinsic diameter is that  in any filling (that is, any \emph{van~Kampen diagram}), there is a family of concentric rings of 2--cells (specifically \emph{$t$--rings}, as defined in Section~\ref{corridors and rings}) that nest to a depth of approximately $n^m$.  
 
\begin{figure}[ht]
\psfrag{D}{$\Delta$}
\psfrag{t}{$t$--arc of}
\psfrag{l}{length $\sim\!n^m$}
\psfrag{n^m}{$\sim\!n^m$}
\psfrag{s}{$\sim n^{m/k}$}
\centerline{\epsfig{file=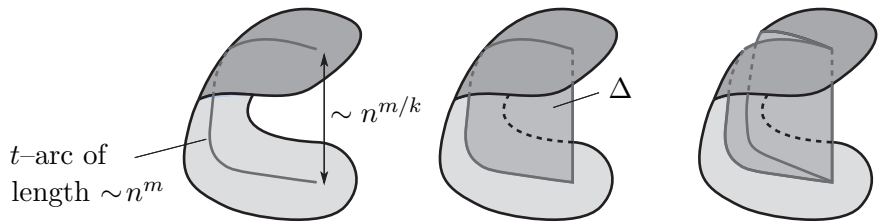}} 
\caption{Inserting a \emph{shortcut diagram} $\Delta$.}
\label{shortcutting fig}
\end{figure}

Let us turn our attention to the role of $\Phi_k$. For each integer $n$ there are ``shortcut words'' of length roughly  $n^{m/k}$ that equal $t^{n^m}$ in $\Phi_k$.   In the Cayley 2--complex of $\Psi_{k,m}$ these mitigate the effect of the nested $t$--rings and cause our large intrinsic--diameter diagrams to have smaller extrinsic diameter.  This is illustrated in the leftmost diagram of Figure~\ref{shortcutting fig}.  

However there remains a significant problem.  In the Cayley 2--complex of $\Phi_k$, the path  corresponding to $t^{n^m}$ and the path of its shortcut word, together form a loop.  This loop can be filled.  We call any filling-disc for this loop a \emph{shortcut diagram} --- $\Delta$ of the middle diagram of Figure~\ref{shortcutting fig} is an illustration.  One could insert two back--to-back copies of the shortcut diagram to get a new filling--disc for our original loop, as illustrated in the rightmost diagram of the figure.  The danger is that this lowers intrinsic diameter and thereby stymies our efforts to separate the two diameter functionals of $\Psi_{k,m}$.  But $\Phi_k$ is built in such a way that  the shortcut diagrams  are  \emph{fat} --- that is, they themselves have large intrinsic diameter (see Section~\ref{are fat}).  So  inserting shortcut diagrams 
decreases  \emph{intrinsic} diameter far less than the presence of shortcut words decreases extrinsic diameter.  

We would like to say that the upshot is that for $\Psi_{k,m}$, the intrinsic diameter functional is at least  $n^m$ and the extrinsic diameter functional at most $n^{m/k}$. But, in truth, technicalities make our bounds, detailed in Theorem~\ref{Main thm with details}, more complicated.

A noteworthy aspect of the construction of $\Phi_k$ is the use 
of an \emph{asymmetric} HNN 
extension $J$ (cf.\ Section~\ref{Phi}).  We anticipate such extensions will prove useful in other contexts.
We build  $\Phi_k \cong \Theta_k \ast_{\Z} J \ast_{\Z} \widehat{\Theta}_k$ by amalgamating $J$ with groups $\Theta_k \cong \widehat{\Theta}_k$ of the form $\Z^k \rtimes F_2$.  
We build  $\Gamma_m$ by starting with a finitely generated free group and taking a number of HNN extensions and amalgamated free products.  All of these groups have compact classifying spaces; $\Psi_{k,m}$ has geometric dimension $k+1$.
The Dehn function of $\Psi_{k,m}$ is at most $n \mapsto C^{n^k}$ for some constant $C>0$ --- see Remark~\ref{main Dehn fn rem}. 
\ms

\subsection{The organisation of this article}

In Section~\ref{Riemannian-algebraic} we prove Theorem~\ref{Riem-Comb},
establishing the qualitative agreement of the group theoretic and Riemannian definitions of diameter.
In Section~\ref{the groups} we present the groups
$\Psi_{k,m} = \Phi_k \ast_\Z  \Gamma_m$ that will be used
to prove Theorem~\ref{main thm}.  Section~\ref{corridors and rings} is a brief discussion of \emph{rings} and \emph{corridors} --- key tools for analysing van~Kampen diagrams.  In Sections~\ref{Phi_k salt} and \ref{d and d2} we explain the salient isodiametric and distortion properties of $\Phi_k$ and $\Gamma_m$.
In Section~\ref{main thm section} we prove the main theorem, modulo generalities about the diagrammatic
behaviour of retracts and amalgams, which we postpone to  Section~\ref{Amalgams and retractions}.  Section~\ref{qi invariance} is dedicated to a proof that, up to $\simeq$ equivalence, extrinsic and intrinsic diameter are quasi--isometry invariants amongst finitely presentable groups.

\section{Riemannian versus combinatorial diameter} \label{Riemannian-algebraic}

In this section we reformulate the combinatorial intrinsic and extrinsic diameter filling functions for finitely presented groups $\Gamma$ in terms of van~Kampen diagrams, and then relate them to their Riemannian analogues.  

\subsection{Diameter of van~Kampen diagrams} \label{van Kampen diameters}
  
Let $\PP= \langle \AA \mid \RR \rangle$  be a finite presentation
of a group $\Gamma$. 
Let $\Delta$ be a van~Kampen diagram over $\PP$ with base vertex $\sstar \in \partial \Delta$.  
(We assume that the reader is familiar with  basic definitions and properties of 
van Kampen diagrams and Cayley 2--complexes --- \cite{Bridson6} is a recent survey.)   Let $\rho$ be the path metric on $\Delta^{(1)}$ in which every
edge has length one and let $d$ denote the word metric on $\Gamma$ associated to $\AA$. 
Let $\mathcal C$ be the Cayley graph of $\Gamma$ with respect to
$\AA$ and let $\pi:\Delta^{(1)}\to\mathcal C$ be the label--preserving
graph--morphism with $\pi(\star)=1$. Define the \emph{intrinsic} and \emph{extrinsic diameter} of $\Delta$ by
\begin{eqnarray*}
\IDiam(\Delta) & := & \max \set{ \rho(\sstar,v) \mid v \in
\Delta^{(0)} },  \\
\EDiam_{\PP}(\Delta) & := & \max \set{ d(1,\,\pi(v)\,) \mid v
\in \Delta^{(0)} }.
\end{eqnarray*}

We shall show  that the algebraic definitions of diameter for
 null--homotopic words given
 in Section~\ref{intro}
are equivalent to
\begin{eqnarray*}
\IDiam_{\PP}(w) & := & \min \set{\IDiam(\Delta) \mid \Delta \text{ a
van~Kampen diagram for } w}, \textup{and} \\
\EDiam_{\PP}(w) & := & \min \set{\EDiam_{\PP}(\Delta) \mid \Delta \text{ a
van~Kampen diagram for } w}.
\end{eqnarray*}
  
Given a van~Kampen diagram $\Delta$ for $w$, one can cut it open
 along  a maximal tree $\TT$ in $\Delta^{(1)}$  to produce a ``lollipop'' diagram 
$\Delta_0$ (with a map $\pi_0: \Delta_0 \to \mathcal C$)
whose boundary circuit is labelled by a word
${w_0} = \prod_{i=1}^N {u_i}^{-1}r_iu_i$ that equals $w$ in $F(\AA)$ and has $r_i \in \RR^{\pm 1}$ for all $i$.  There is a natural combinatorial folding map $\varphi : \Delta_0 \onto \Delta$ such that $ \pi \circ \varphi = \pi_0$.  

It follows that
$\max \set{ \, d(1, \bar{p}) \, \mid \, p \textup{ a prefix of } w_0 \, }$, where $\bar{p}$ denotes the element of $\Gamma$ represented by $p$, is at most $\EDiam_{\PP}(\Delta)$, because $p$ corresponds to a vertex on $\partial \Delta_0$.  
Thus the algebraic version of $\EDiam_{\PP}(w)$ is bounded above by the geometric
version.  For the opposite inequality,  suppose that $x$ is a word of the form $\prod_{i=1}^N {u_i}^{-1} r_i u_i$ where $r_i \in \RR^{\pm 1}$ for all $i$, and $x=w$ in $F(\AA)$.  Define  $m:= \max \set{ \, d(1, \bar{p}) \, \mid \, p \textup{ a prefix of } x \, }$.  Suppose $x$ minimises $m$ and is of minimal length amongst all words that minimise $m$.  
There is an obvious
lollipop diagram with boundary word $x$. A van~Kampen diagram $\Delta$ for $w$ with $\EDiam_{\PP}(\Delta) \leq m$ can be obtained from this by 
successively folding together pairs of edges with common initial
point and identical labels.  (A concern is that in the course of this folding, 2--spheres might be pinched off, but minimising $\ell(x)$ avoids this.)  

\medskip

$\IDiam_{\PP}(w)$, as defined above, differs from $\IDiam_{\PP}(w)$ of Section~\ref{intro} by at most   $\frac 1 2\max\set{\ell(r) \mid r \in \RR} + \ell(w)$. This can be shown by an argument similar to the one above, expect that for the first inequality the tree  $\TT$ should be
chosen to be a maximal geodesic tree based at $\sstar$, and when
proving the reverse inequality one must specify first that $x$ minimises 
$\max_i \ell(u_i)$ and then that it has minimal length subject to this
constraint.  This discrepancy is of no great consequence for us since it has no effect on the $\simeq$--class of the resulting functional.  

\medskip
      
We will use the van~Kampen--diagram interpretations of $\IDiam_{\PP}(w)$ and $\EDiam_{\PP}(w)$ in the remainder of this article.  
  
\subsection{Generality}  \label{regularity subsection}

With an eye to future applications and to highlight the essential ingredients of our arguments, in this subsection and the next  we will work in a more general setting than that of Riemannian manifolds.

\begin{defn}
Define the intrinsic and extrinsic diameters of a disc--filling $F:\mathbb{D}^2\to X$ in a metric space $(X,d_X)$ by \begin{eqnarray*}
\IDiam(F) \!\! &  := & \!\! \sup \set{ d_F(a,b) \ \mid \ a,b \in \mathbb{D}^2 \, }, \text{ and} \\ 
\EDiam(F) \!\! &  := & \!\! \sup \set{ d_{X}(a,b) \ \mid \ a,b \in \mathbb{D}^2 \, }
\end{eqnarray*}
where $d_F$ is the pseudometric $$d_F(a,b) := \inf \set{\ell(F \circ p) \mid p \textup{ a path in } \mathbb{D}^2 \textup{ from } a \textup{ to } b }$$ on $\mathbb{D}^2$.  Define
$\IDiam_{X}, \EDiam_{X} : [0,\infty) \to[0,\infty)$ by
\begin{eqnarray*}
\IDiam_{X}(l) \!\! & := & \!\! \sup_c \inf_{F} \left\{ \ \IDiam(F) \left| \textup{ disc--fillings } F \textup{ of loops } c \textup{ with } \ell(c) \le l  \ \right. \right\}, \\
\EDiam_{X}(l) \!\! & := & \!\! \sup_c \inf_{F}  \left\{ \ \EDiam(F) \left| \textup{ disc--fillings } F \textup{ of loops } c  \textup{ with } \ell(c)\le l \ \right. \right\}.
\end{eqnarray*}
\end{defn}

[In closer analogy with the definitions
in the previous subsection, one could deal instead with based loops and discs
and define diameter functions in terms of distance from the basepoint. This
makes little difference: the resulting  functions $[0,\infty)\to
[0,\infty)$ are $\simeq$.]

\begin{lemma}   \label{well defined} Let 
 $X$ be the universal cover of a compact geodesic space $Y$ for which  there exist $\mu, A > 0$ such that every loop of length less than $\mu$ admits a disc--filling of intrinsic diameter less than $A$. Equip $X$
with the induced length metric.
 Then $\IDiam_X$ and $\EDiam_X$ are well--defined functions 
 $[0,\infty)\to [0,\infty)$.
\end{lemma}

\Proof
First observe that every rectifiable loop in $X$ admits a disc--filling of finite intrinsic diameter: given a disc--filling  $F:\mathbb D^2\to X$ of a rectifiable loop one can triangulate the disc so that the image of each edge has diameter at most $\mu/3$, and then one can modify $F$ away from $\partial \mathbb{D}^2$ and the vertex set so that its restriction to each internal edge is a geodesic and every triangle has intrinsic diameter at most $A$.
 
Next, suppose $c$ is a loop in $X$ of length $l$.
Shrinking $\mu$ if necessary, we may assume that balls of radius $\mu$ in $Y$ lift to $X$.
Cover $Y$ with a maximal collection of disjoint balls of radius $\mu/10>0$;  let $\Lambda \subset X$ be the set of lifts of their centres.    Partition $c$ into $m \leq 1+10\ell/\mu$ arcs, each of length at most $\mu/10$ and with end points $v_0, v_1, \ldots, v_m$.  Each $v_i$ is a distance at most  $\mu/10$ from a point $u_i \in \Lambda$, and the distance from $u_i$ to $u_{i+1}$ is less than $\mu/2$ (indices modulo $m+1$).  The loops made by joining $v_i$ to $u_i$, then $u_i$ to $u_{i+1}$,  then $u_{i+1}$ to $v_{i+1}$ (each by a geodesic), and then  $v_{i+1}$ to $v_i$ by an arc of $c$, each have total length at most $\mu$.  So, by hypothesis, they admit  disc--fillings with intrinsic diameter at most $A$.  Such disc--fillings together form a collar between $c$ and a piecewise geodesic loop formed by concatenating at most $1+10\ell/\mu$ geodesic segments, each of length at most $\mu/2$ with endpoints in $\Lambda$.   Modulo the action of $\pi_1Y$,  there are only finitely many such piecewise geodesic loops, and by the argument in the previous paragraph each one admits a filling of finite intrinsic diameter.  It follows that  
$\IDiam_X$ and (hence) $\EDiam_X$ are well--defined functions. 
\qed

\begin{remark} \label{M good} 
If $X= \widetilde{M}$, the universal cover of a closed connected Riemannian manifold $M$, then  it satisfies the conditions of Lemma~\ref{well defined}.
Indeed, the required $\mu$ and $A$ exist for any cocompact space $X$ that is locally uniquely geodesic, for example a space with upper curvature bound in the sense of A.D.~Alexandrov --- i.e.\ a
CAT$(\kappa)$ space \textup{\cite{BrH}}: by cocompactness, there exists $\eta >0$
such that geodesics are unique in balls of radius $\eta$, and any loop
in such a ball can be filled  by coning
it off to the centre of the ball using geodesics.
\end{remark}

\subsection{The Translation Theorem }\label{s:reduce} \label{ss23} 

Theorem \ref{Riem-Comb} 
is a special case of the following result.

\begin{thm} \label{t:approx}
Suppose a group $\Gamma$ with finite presentation $\PP$ acts properly and cocompactly by
isometries on a simply connected geodesic metric space $X$ for which there exist $\mu, A > 0$ such that every loop of length less than $\mu$ admits a disc--filling of intrinsic diameter less than $A$. 
Then $\rm{IDiam}_{\PP} \simeq \rm{IDiam}_X$ and $\rm{EDiam}_{\PP} \simeq \rm{EDiam}_X$.
\end{thm}

We shall first establish the assertion concerning extrinsic diameter and the relation $\rm{IDiam}_{X} \preceq \rm{IDiam}_{\PP}$.

Map the Cayley graph of $\G$ to $X$ as follows:  fix a basepoint $p \in X$, choose a geodesic from $p$ to its translate $a\cdot p$  for each generator $a$
of $\PP$, and then
extend equivariantly. Following \cite{Bridson6},  a path in $X$ is called {\em word--like} if it is the image in $X$ of an edge--path in the Cayley graph.

Each 2--cell in the Cayley 2--complex $\Cay^2(\PP)$ is attached to the Cayley graph by an edge--loop
labelled by one of the defining relations $r$ of $\PP$. For each $r$
we choose a filling disc $F_r$ of
finite intrinsic diameter for the corresponding word--like loop in $X$ based at $p$. 
We then map
$\Cay^2(\PP)$ to $X$ by the $\G$--equivariant map that sends each
2--cell with boundary label $r$ to a translate of $F_r$.

Using a collar between an arbitrary rectifiable loop in $X$ and a word--like loop, as in the proof of Lemma~\ref{well defined}, one can show that there is no change in the $\simeq$ classes of either $\rm{IDiam}_X$ or
$\rm{EDiam}_X$ if one takes the infima in their definitions to be over fillings of
word--like loops only.  
Having made this reduction, the relation $\rm{IDiam}_X  \preceq \rm{IDiam}_{\PP}$ becomes obvious, since one gets an upper bound
on the intrinsic diameter of discs filling a word--like loop
simply by taking  the image in $X$ of a minimal
intrinsic diameter van~Kampen diagram $\Delta$ for the appropriate word.  A technical concern here is that $\Delta$ need not be a topological 2--disc, but rather may be a finite planar tree--like arrangement of topological 2--discs and 1--dimensional arcs. This can be overcome by extending $\Delta\to
X$ to a small regular neighbourhood $V$ of $\Delta\subset\mathbb R^2$ by
a map  that is constant on the slices $\{x\}\times [0,1)$ of the annulus
$V\smallsetminus \Delta \approx \partial V\times [0,1)$.

The same argument yields $\rm{EDiam}_X
 \preceq \rm{EDiam}_{\PP}$, and the converse is an
 easy approximation argument (the details of which are
 included in Lemma \ref{l:approx}). Briefly, noting that
 the map of the Cayley graph of $\G$ to $X$ is a quasi--isometry,
 we need only show that  a disc--filling $F:\mathbb D^2\to X$ of a 
 word--like loop gives rise to a van Kampen diagram for the corresponding
 word that is $C_0$--close to $F$. Such an approximating
 disc is obtained by simply taking a fine triangulation of
 $\mathbb{D}^2$, and then labelling it as in the proof of Lemma~\ref{enough}.
 
\bs

The remainder of this section is devoted to
the most difficult relation in Theorem~\ref{t:approx}, namely
$\rm{IDiam}_{\PP} \preceq \rm{IDiam}_X$. Like the argument sketched above, the
 proof of this assertion involves constructing a suitable combinatorial
 approximation to a Riemannian filling of a word--like loop. The control
 required in this approximation is spelt out in the following lemma.

\begin{lemma} \label{enough} \label{l:approx} 
To prove $\IDiam_{\PP} \preceq \IDiam_X$, it suffices to exhibit  constants $L_0,L_1,\lambda$ such that, given a  filling $F:\mathbb{D}^2\to X$ of intrinsic diameter $\delta$ for a word--like loop $c$, one can construct a  combinatorial cellulation $\Delta$ of $\mathbb{D}^2$ and a map $f: \Delta^{(0)} \to X$ with the following properties:
\begin{enumerate}
\item[\textup{(}\textit{1}\textup{)}] the attaching map of each 2--cell of $\Delta$ has combinatorial length at most $L_0$,
\item[\textup{(}\textit{2}\textup{)}] adjacent vertices in $\Delta$ are mapped by $f$ to points that are a distance
at most $L_1$ apart in $X$,
\item[\textup{(}\textit{3}\textup{)}] from each vertex of $\Delta$, one can reach 
 $\partial \Delta$ by traversing a path consisting of at most $\lambda (\delta+1)$
 edges, and
 \item[\textup{(}\textit{4}\textup{)}] the vertices $\partial\Delta^{(0)}$ on $\partial \Delta$ are mapped by $f$ to points on $c$, and their cyclic ordering is
preserved; moreover, every vertex of $c$ is in $f(\partial\Delta^{(0)})$.
 \end{enumerate}
\end{lemma}

\Proof In light of Theorem \ref{qi invariance thm}, we may
take $\PP$ to be a finite presentation well adapted to the geometry of $X$.
Thus we fix $\xi >0$ so that $X$ is the $\xi$--neighbourhood of  $\Gamma \cdot p$,
and as generators of $\G$ we take the set $\mathcal A$
 of those $\gamma$ such that
$d(p,\, \gamma(p)) \le L_1 + 2\xi$. As relations $\mathcal R$
we take all words of 
length at most
 $L_0$ in the letters $\mathcal A$ that equal the identity in $\Gamma$.
The remainder of our proof shows that these relations suffice;
cf.~\cite[page 135]{BrH}.

Given a word $w$ in the generators $\mathcal A$ that equals $1$ in $\G$, we
consider the corresponding word--like loop $c$ in $X$. By construction,
there are constants $\alpha, \beta$ such that 
the length of $c$ is at most $\alpha \ell(w) + \beta$.  So we will be
done if, given a disc--filling of $c$ with intrinsic diameter $\delta$,
we can exhibit a van Kampen diagram for $w$ with intrinsic diameter
bounded above by a linear function of $\delta$. 

Let $\Delta$ be as described in the statement of the
lemma and label
 each vertex $v\in\Delta^{(0)} \ssm \partial \Delta$ by a group element $\gamma_v$
 such that $\gamma_v(p)\in \G \cdot p$ is within
 $\xi$ of $f(v)$. Label an
 edge from a vertex $v$ to a vertex $v'$ by $\gamma_v^{-1}\gamma_{v'}$, which is in $\mathcal A$. For vertices $v$
 of $\partial \Delta$, we choose $\gamma_v$ to
 be one of the endpoints of the edge of the word--like loop
 $c$ in which (the image of) $v$ lies.
 
 By hypothesis,
 the boundary of each 2--cell in the cellulation is then labelled
 by a word that belongs to $\mathcal R$. This process
 does not quite yield a van Kampen diagram for the original word $w$, but
 rather for a word $w'$ that is equal to $w$ modulo the
 cancellation of edge--labels $1$. (Our convention on the choosing of $\gamma_v$
 for vertices on $\partial\mathbb D^2$ ensures this equality.)
 Collapsing such edges yields the
 desired diagram.
 \qed
\bs

Continuing the proof that  $\IDiam_\PP\preceq \IDiam_X$, we concentrate on a disc--filling $F:\mathbb{D}^2 \to X$ of intrinsic diameter $\delta$ for a word--like loop $c$ with basepoint $F(\star)=p$.   We will construct $f$ and $\Delta$ satisfying the conditions of Lemma~\ref{l:approx} with $L_0 =  6$, $L_1 = 5/2$ and $\lambda = 1$ in three steps.

\ms 

\ni (i) \emph{Constructing the tree of a fine triangulation.}

\ss

\def\e{\varepsilon}

Since $F$ is continuous and $\mathbb{D}^2$  is compact, $\mathbb{D}^2$ admits a finite triangulation $\hat{\Delta}$ such that $F$ maps each triangle to a subset of $X$ of diameter at most $1/2$.  Fix $\varepsilon >0$.  Fix an ordering $v_0,v_1,v_2,\dots$ on the vertices $V$ of $\hat{\Delta} \ssm \partial \hat{\Delta}$ and, proceeding along the list, choose an embedded arc $\sigma_v$ in $\mathbb{D}^2$ joining $v$ to $\star$ as follows so that $T := \bigcup_{v \in V} \sigma_{v}$ is a tree.  Having defined $\sigma_{v_i}$ for each $i<m$, to define $\sigma_{v_m}$ first choose an embedded arc from $v_m$ to $\star$ of length at most $\varepsilon + d_F(\star, v_m)$, and 
arrange that it pass through no vertices in $V \ssm \set{ v_m , \star}$
(more on this in a moment) and finally replace its terminal segment beginning at the first point $x$ it meets in $\bigcup_{i<m} \sigma_{v_i}$, with the arc in  
$\bigcup_{i <m} \sigma_{v_i}$ between $x$ and $\star$.  By construction $V$ is the set of leaves of $T$.  Were it not for the $\varepsilon$--error terms, the diameter of $T$ would be at most $2 \delta$. We choose $\varepsilon$ small enough so that the diameter is no more than $2 \delta +1$.

[Technicalities: We introduced $\varepsilon$ to avoid a discussion of the
existence of geodesics in the pseudometric space $(\mathbb D^2, d_{F})$.
We can arrange that $\sigma_{v_m}$ not pass through $V \ssm \set{ v_m , \star}$
by replacing $F$ with a map $F'$ that is constant on 
a small disc about each $v\in V$ and which, on the complement $E$
of  the union of these disjoint
discs, is $F\circ h$ where $h: E\to \mathbb D^2$ 
is a homeomorphism stretching an annular neighbourhood of each deleted disc
to a punctured disc with the same centre. The distances between
distinct vertices in  
$(\mathbb D^2, d_{F})$ and $(\mathbb D^2, d_{F'})$ are the same
and, in the latter, a path can take a detour near any vertex without adding
any length.]

We intend to combine $T$ and $\hat{\Delta}$ to give a cellulation $\Delta$ with the properties required for Lemma~\ref{enough}.  However, as things stand, the intersection of $T$ with the cells of $\hat{\Delta}$ could be extremely complicated.  To cope with this we introduce a hyperbolic device to extract the essential features of  the pattern of intersections.

\ms

\ni (ii) \emph{Taming the intersection of $T$ with $\hat{\Delta}$.}

\ss

As above, $V$ is the set of vertices of $\hat{\Delta}$ that lie in the interior of $\mathbb{D}^2$.
We impose a complete hyperbolic metric on $\mathbb{D}^2 \ssm V$ in such a way that  $\partial \mathbb{D}^2$ is a regular geodesic polygon and the 1--cells of $\hat{\Delta}$ are geodesics in the hyperbolic metric.  For each $v \in V$ let $\tau_v$ be the unique geodesic in the same homotopy class of ideal paths from $v$ to $\star$ in $\mathbb{D}^2 \ssm V$ as $\sigma_v$.
To aid visualisation, we are going to talk of these $\tau_v$ arcs as being \emph{red}.

\begin{figure}[ht]
\psfrag{v}{$v_m$}
\centerline{\epsfig{file=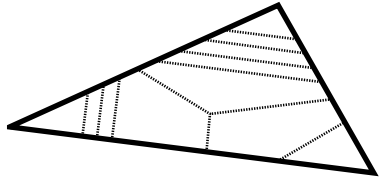}} \caption{The intersections of $\bigcup_{v \in V} \tau_v$ with a triangle of $\hat{\Delta}$.} \label{good triangle}
\end{figure}

The point of what we have just done is that it vastly simplifies the pattern of intersections of
the $\sigma_v$ with $\hat{\Delta}^{(1)}$: the pattern of intersections of the $\tau_v$ with a triangle is the simple pattern shown in Figure~\ref{good triangle}  (except if the triangle has $\star$ as one of its vertices then some of the red edges can meet $\star$, indeed the sides of the triangle may be red).
Moreover, the larger--scale geometry has been retained in a way that allows us to construct a cellulation with combinatorial diameter at most a linear function of $\delta$, as we will explain. 

\ms

\ni (iii) \emph{Constructing the good cellulation $\Delta$ and map $f$.}

\ss

\def\t{\tau}

We begin by taking the cellulation $\Delta_0$ consisting of $\hat{\Delta}$ overlaid with $\hat{T} := \bigcup_{v \in V} \tau_v$.  That is, the vertices are the vertices of $\hat{\Delta}$ and the points of intersection of $\hat{T}$ with edges of $\hat{\Delta}$, and the edges come in two types:  \emph{black} edges which are edges of $\hat{\Delta}$ (which may be partitioned into multiple edges by vertices on $\hat{T}$), and \emph{red} edges, which are arcs in $\hat{T}$ between vertices.   Define $f_0 : {\Delta_0}^{(0)} \to X$ to be the restriction of $F$ to ${\Delta_0}^{(0)}$.  

We will alter  $\Delta_0$ (and $f_0$) and the paths $\tau_v$ until we can control their combinatorial lengths in terms of the lengths of the $\sigma_v$ and thereby get an upper bound on the  intrinsic diameter of $\Delta_0$.

\begin{figure}[ht]
\psfrag{v}{$v$}
\psfrag{w}{$v'$}
\psfrag{s}{$\star$}
\centerline{\epsfig{file=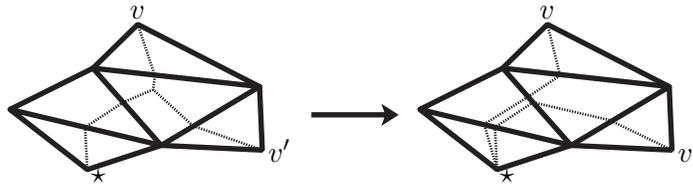}} \caption{Doubling along the intersection of two paths $\tau_v$ and $\tau_{v'}$.} \label{double fig}
\end{figure}

It will be convenient if whenever $v, v' \in V$ are distinct, $\tau_v$ and $\tau_{v'}$  meet only at $\star$.  Achieve this by doubling $(\tau_v \cap \tau_{v'}) \ssm \set{\star}$, and joining doubled vertices by a new black edge --- see Figure~\ref{double fig}. Call the new combinatorial disc $\Delta_1$ and define $f_1 : {\Delta_1}^{(0)} \to \widetilde{M}$ by letting $f_1(u)$ be $f_0(\bar{u})$ where $\bar{u}$ is the image of $u$ under the obvious combinatorial retraction of $\Delta_1$ onto $\Delta_0$.

Repeat the following for each $v \in V$ such that $\tau_v$ has combinatorial length at least 2. 
As $\tau_v$ is in the same homotopy class of ideal paths from $v$ to $\star$  of $\hat{\Delta} \ssm V$ as $\sigma_v$,  the sequence of edges $\tau_v$ crosses is also crossed by $\sigma_v$ in the same order, except $\sigma_v$ may make additional intersections en route.  That is, if $\set{x_i}_{i=0}^k$ are the points of intersection of $\tau_v$ with  $\hat{\Delta}$  numbered in order along $\tau_v$ with $x_0 = \star$, $x_k=v$ and all other $x_i$ on
so a black edge $e_i$, then there are points $\set{y_i}_{i=0}^k$ occurring in  order along $\sigma_v$ with $y_0 = x_0$, $y_k = x_k$ and $y_i$ on $e_i$ for all other $i$.   For all $i$ let $n_i$ be the integer part of the length of the segment of $\sigma_v$ from $y_i$ to $\star$.  Note that $i \mapsto n_i$ increases monotonically from $0$ to $\lfloor\ell(\sigma_v)\rfloor$, and takes all possible integer values in that range because of the condition that the diameters of the images of the triangles in $\hat{\Delta}$ are at most $1/2$.
For each $n \in \N$, collapse to a single vertex the minimal  arc of  $\tau_v$ that includes all red vertices $x_i$ such that $n_i = n$.  Note that as  each arc is collapsed the complex remains topologically a disc since the arc does not run between two boundary vertices.   

The result is the desired cellulation $\Delta$ of $\mathbb{D}^2$ and we define $f : {\Delta}^{(0)} \to \widetilde{M}$ by setting $f(\bar{u})$ to be $f(u)$ for any choice of preimage $u$ of $\bar{u}$ under the quotient map $\Delta_1 \onto \Delta$.   All that remains is to explain why $\Delta$ and $f$ satisfy the conditions of Lemma~\ref{enough} with $L_0 =6$, $L_1 = 5/2$, and $\lambda = 1$.

The faces of $\Delta$ each have at most six sides since this was the case for $\Delta_1$ (six being realised when red edges cut across each of the corners of a triangle of $\hat{\Delta}$) and the subsequent collapses of edges can only decrease the number of sides, which proves (\emph{1}). 

For (\emph{2}), first suppose $x$ and $x'$ are vertices in $\Delta$ joined by a red edge.  Then they are images of vertices $x_i$ and $x_{i'}$ on the some path $\tau_v$ in $\Delta_1$ such that $f(x) = f_1(x_i)$, $f(x') = f_1(x_{i'})$ and $\abs{n_i - n_{i'}} =1$.  But then, as $\abs{n_i - n_{i'}} =1$,  
$$d(f(x), f(x'))  \ = \  d(f_1(x_i) , f_1(x_{i'})) \ \leq \  2.$$  
Next suppose $x$ and $x'$ are vertices in $\Delta$ joined by a black edge.  Assume they are images of vertices $x_i$ and $x_{i'}$ on paths $\tau_v$ and $\tau_{v'}$ in $\Delta_1$ such that $x_i$ and $x_{i'}$ are joined by a black edge.  Then $d(f_1(x_i) , f_1(x_{i'})) \leq 1/2$ and there are $x_j$ and $x_{j'}$ on $\tau_v$ and $\tau_{v'}$ such that $f(x) = f_1(x_j)$, $f(x') = f_1(x_{j'})$, $n_i= n_j$, and $n_{i'} = n_{j'}$.  So 
\begin{eqnarray*}
d(f(x), f(x'))  & = &  d(f_1(x_j) , f_1(x_{j'})) \\  & = & d(f_1(x_j) , f_1(x_{i})) + d(f_1(x_{i}) , f_1(x_{i'})) + d(f_1(x_{i}') , f_1(x_{j'})) \\ & \leq &  1 + \frac{1}{2}  + 1 \ = \ \frac{5}{2}.
\end{eqnarray*}
The remaining possibility is that $x$ and $x'$ are joined by a black edge but one or both is in $\partial \Delta \ssm \set{\star}$ and so is not on any $\tau_v$ in $\Delta_1$. Similar considerations to the previous case yield $d(f(x), f(x')) \leq 1 + 1/2$.   

By construction, for all $v \in V$, the combinatorial length of the images in $\Delta$ of $\tau_v$ is at most the length of $\sigma_v$, and so at most $\delta+1$ and we have  (\emph{3}).  And finally, (\emph{4}) is immediate from the way we constructed $f$ and $\Delta$.  

\section{The groups for Theorem~\ref{main thm}}\label{the groups}

The groups that we use form a
2--parameter family
$$\Psi_{k,m} \ = \  \Phi_k\ast_\Z\Gamma_m,$$
where  $m>k>1$ are integers.

\subsection{The ``Fat--Shortcuts'' Groups $\Phi_k$}\label{Phi}

We construct $\Phi_k$ by amalgamating $J$, a double HNN--extension of $\Z$,
with two isomorphic copies of 
$\Theta_k=\Z^k\rtimes F_2$, where the action of $F_2$ is chosen
so that a certain cyclic subgroup $\langle s_k \rangle$ of $\Z^k$
is distorted in a precisely controlled manner as will be required in Lemma~\ref{s_k distortion}.   
The two copies  of $\Theta_k$ are amalgamated with
 $J$ along the cyclic subgroups $\langle s_k \rangle$. [The distortion of a finitely generated
subgroup $S$ of a finitely generated group $G$ measures the difference
between the word metric on $S$ and the restriction to $S$ of the word
metric of $G$ --- see \cite{BrH}, page~506.]

\ms

We first take a  presentation $\OO_k$ of
 $\Theta_k = \Z^k\rtimes F_2$ by choosing a basis
$s_1, \ldots, s_k$ for $\Z^k$ and specifying the action of $F_2=\langle f,g\rangle$:
\begin{tabbing}
A \= AAAAAAAA \= AAAAAAAAAAAAAAAAAAAAAAAA  \kill %
\> \emph{generators} \> $s_1,  \ldots, s_k, f, g$ \\ %
\> \> \\ %
\> \emph{relations} \> $f^{-1}s_k f=s_k, \ \textup{and} \ \forall\, i < k, \ f^{-1} s_i f=s_i s_{i+1} $ \\  %
\> \> $g^{-1}s_k g=s_k, \ g^{-1}s_{k-1} g=s_{k-1}, \ \textup{and} 
\ \forall\, i < k-1, \ g^{-1} s_i g=s_i s_{i+1}  $ \\  %
\> \> $\forall\, i \neq j, \ [s_i, s_j]=1 \,.$ %
\end{tabbing}

\ni Let $\widehat\Theta_k$ be a second copy of $\Theta_k$
with corresponding presentation $\widehat{\OO}_k$ in which the generators are $\hat{s}_1,  \ldots, \hat{s}_k, \hat{f}, \hat{g}$.

\ms We define $J$ to be the amalgam along $\langle b \rangle \cong \Z$ of 
the {\em asymmetric  HNN extension}  $\langle b,t, s \mid
{t}^{-1} b s = b^3, \ {s}^{-1}b t = b^3 \rangle$ with the standard HNN extension $\langle b, \hat s \mid 
\hat{s}^{-1} b \hat{s}= b^3 \rangle$.  [This is obtained from
$\langle b,t\mid t^{-1}b^2t=b^6\rangle$ by the
Tietze move introducing $s=b^{-1}tb^3$.]  Thus $J$ has presentation
$$\langle b,t, s,\hat s \mid {t}^{-1} b s = b^3, \ {s}^{-1}b t = b^3,  
\hat{s}^{-1} b \hat{s}= b^3 \rangle.$$

 \ni Finally, we define
$$\Phi_k \ := \ \Theta_k \ast_\Z J \ast_\Z\widehat{\Theta}_k,$$ 
where the amalgamations identify  $s,\hat{s}\in J$ with $s_k\in\Theta_k$
and $\hat s_k\in\widehat\Theta_k$, respectively. Thus, writing 
$\OO_k = \langle \Sigma_k \mid \Upsilon_k \rangle$ and $\widehat{\OO}_k = \langle \widehat{\Sigma}_k \mid \widehat{\Upsilon}_k \rangle$, a presentation for $\Phi_k$ is
$$\PP_k \ := \ \langle b,t, \Sigma_k, \widehat{\Sigma}_k\mid \Upsilon_k, \widehat{\Upsilon}_k,\, 
{t}^{-1} b s_k = b^3, \ {s_k}^{-1}b t = b^3, \
\mbox{$\hat{s}_k$}^{-1} b \hat{s}_k = b^3 \rangle.$$

The asymmetric nature
of the quasi--HNN presentation of $\langle b,s,t\rangle\subset J$ will be important
in Section~\ref{shortcutting} when we wish to force shortcuts to be fat  --- see  Figure~\ref{Delta_m} for an example and Section~\ref{shortcuts are fat} for the general argument. As we indicated in the introduction, this fatness
 is the key to forcing the behaviour of the intrinsic and extrinsic diameter functionals to diverge.  

\subsection{The $t$--rings groups $\Gamma_m$} \label{Gamma}

\def\s{\sigma}

The groups  $\Gamma_m$, whose presentations $\QQ_m$ we give below, are 2--fold HNN
extensions of the free--by--free groups $B_m=F_m\rtimes F_2$. By
deleting from  $\QQ_m$  the generator $T$ and the defining relations in which it appears,  
one recovers the 
groups constructed by the first author in \cite{Bridson}; these
have isodiametric properties that prevent their
asymptotic cones from being simply connected.

We define $B_m:=F_m\rtimes F_2$, where the first generator $\s$
of $F_2 = \langle \s, t \rangle$ acts on $F_m$ as an automorphism with polynomial
growth of degree $m-1$, and $t$ acts trivially. Specifically, $B_m$ has presentation
$$
\mathcal B_m \ := \  \langle \AA_m, \s, t \mid \s^{-1} a_m \s= a_m\,; \ \forall\, i <m, \ 
\s^{-1} a_i \s = a_i a_{i+1}; \forall\, j, \ [t,a_j] =1\rangle,
$$
where $\AA_m:= \set{a_1, \ldots, a_m}$.
We then obtain
$\Gamma_m$ by taking two successive HNN extensions of $B_m$:  in the first, the stable 
letter $\tau$ commutes with a skewed copy of $F_m$ and, in the second, the stable letter $T$ commutes with $\tau$ and $t$.  Thus, abbreviating
$\mathcal B_m = \langle \AA_m, \s, t \mid \RR\rangle$, we define $\G_m$
to be the group with (aspherical) presentation
$$\QQ_m  \ :=  \  \langle \AA_m, \s, t, \tau, T \mid \RR, [t,T]=1, [\tau, T]=1, [\tau, a_m t]=1; \ \forall\, i<m, \ [\tau,a_i]=1\rangle.
$$

\subsection{Preliminary lemmas}

The following two lemmas are self--evident.

\begin{lemma} $\Gamma_m$ retracts onto $B_m=F_m\rtimes\langle \sigma,t\rangle$
and hence onto $\langle t\rangle$.
\end{lemma}

\begin{lemma} \label{two retracts}
The groups $\Theta_k$ and $\hat{\Theta}_k$ 
are retracts of $\Phi_k$ via
\begin{eqnarray*}
\Theta_k \ \stackrel{\iota}{\hookrightarrow} \ \Phi_k &
\stackrel{\psi}{\twoheadrightarrow} & \Theta_k
\\
\hat{\Theta}_k \ \stackrel{\hat{\iota}}{\hookrightarrow} \ \Phi_k
& \stackrel{\widehat{\psi}}{\twoheadrightarrow} & \widehat{\Theta}_k,
\end{eqnarray*}
where $\iota$ and $\hat{\iota}$ are the obvious inclusions,
$\psi$ kills $b, \hat{f}, \hat{g}, \hat{s}_1, \ldots, \hat{s}_k$,
maps $t$ to $s_k$, and $f, g, s_1, \ldots, s_k$ to themselves, and
$\hat{\psi}$ kills $b, t, f, g, s_1, \ldots, s_k$ and maps
$\hat{f}, \hat{g}, \hat{s}_1, \ldots, \hat{s}_k$ to themselves.
\end{lemma}

Since $\OO_k$ is the obvious presentation of the semidirect
product $\Theta_k=\Z^k\rtimes F_2$ we have:

\begin{lemma} \label{abelian subpresentation}
If $w \in \set{ {s_1}^{\pm 1}, \ldots, {s_k}^{\pm 1}}^{\sstar}$ is
null--homotopic in $\OO_k$, then it is also null--homotopic in
$$\langle s_1, \ldots, s_k \mid \forall i \neq j, \ [s_i, s_j]
\rangle \ \cong \ \Z^k.$$
\end{lemma}

Since $\langle b \rangle$ and $\langle 
s_1, \ldots, s_{k-1} \rangle\cong\Z^{k-1}$ do not intersect
the amalgamated subgroups of $\Phi_{k,m}=\Theta_k\ast_\Z J\ast_\Z
\widehat\Theta_k$, they generate their free product. In more
detail:

\begin{lemma} \label{null--homotopic lemma}
If $w \in \set{ b^{\pm 1}, {s_1}^{\pm 1}, \ldots, {s_{k-1}}^{\pm
1}}^{\sstar}$ is a null--homotopic word in $\PP_k$, then it is also
null--homotopic in the subpresentation $$\langle b,
s_1, \ldots, s_{k-1} \mid \forall i \neq j, \ [s_i, s_j] \rangle \
\cong \ \Z \ast \Z^{k-1}.$$
\end{lemma}
 
\begin{lemma} \label{binomial coefficients lemma}
For all $n \in \N$, there exists a \emph{(}positive\emph{)} word $w_n \in \set{s_1, \ldots,
s_k}^{\sstar}$ such that  $w_n=f^{-n} {s_1}^n f^n$ in 
$$ F_k \rtimes \Z \ \cong \ \langle \, s_1, \ldots, s_k, f \mid f^{-1} s_k f = s_k\,; \ \forall i <k, \,
f^{-1} s_i f = s_i s_{i+1} \, \rangle$$  
and
$$w_n \ = \ \prod_{i=1}^{k} {s_i}^{n \TB{n}{i-1}} \ \ \ \text{ in
} \ \ \ \Z^k \ = \ \langle s_1, \ldots, s_k \mid \forall i \neq j, \ [s_i, s_j]  \rangle.$$
\end{lemma}

\Proof Let $w_n$ be the unique positive word in the $s_i$ such that 
$w_n=f^{-n} {s_1}^n f^n$ in $F_k \rtimes \Z$.
The stable letter $f$ acts on the abelianisation $\Z^k$ of $F_k$ as left multiplication by the $k \times k$ matrix with ones on the diagonal
and subdiagonal and zeros elsewhere.  So the result follows from the calculation
$$\left(%
 \begin{array}{cccccc}
  1 &  &  &  &  \\
  1 & 1 &  &  &  \\
   & 1 & \ddots &  &  \\
   &  &   \ddots & 1  &  \\
   &  &   & 1 & 1 \\
\end{array}%
\right) ^{\mbox{$n$}} \ = \  \left(%
\begin{array}{ccccc}
  1 &  &  &  &  \\
  \TB{n}{1} & 1 &  &  &  \\
  \TB{n}{2} & \TB{n}{1} & \ddots &  &  \\
  \vdots & \vdots  &  \ddots & 1  &  \\
  \TB{n}{k-1} & \TB{n}{k-2} &  \cdots  & \TB{n}{1} & 1 \\
\end{array} %
\right).$$
\qed

\section{$t$--corridors and $t$--rings in van Kampen diagrams} \label{corridors and rings}

The use of $t$--\emph{corridors} and $t$--\emph{rings} 
for analysing van~Kampen diagrams is well established. Suppose $\PP=\langle \AA \mid
\RR \rangle$ is a presentation, $\Delta$ is a van Kampen diagram over
$\PP$ and
$t \in \AA$. An edge in $\Delta$ is called a \emph{$t$--edge} if it is labelled by $t$.  Let $\Delta^{\sstar}$ denote the dual graph to the
1--skeleton of $\Delta$ (including a vertex $v_{\infty}$ dual to the region exterior to $\Delta$). 

Suppose $\sigma$ is either a simple edge--loop
in $\Delta^{\sstar}$ and that the edges of $\sigma$ are all
dual to $t$--edges in $\Delta$. Then the subdiagram $C$ of $\Delta$
consisting of all the (closed) 1--cells and 2--cells of $\Delta$ dual to vertices and edges of
$\sigma \ssm \set{v_{\infty}}$, is called a \emph{$t$--ring}  if $\sigma$ does not include $v_{\infty}$
and a \emph{$t$--corridor} in it does.  Let $e_1,
\ldots, e_r$ be the duals of the edges of $\sigma$, in the order
they are crossed by $\sigma$, and such that in the ring case
$e_1=e_r$ and in the corridor case $e_1$ and $e_r$ are in
$\partial e_{\infty}$.  The \emph{length} of the $t$--ring or $t$--corridor is $r-1$ (which is zero for $t$--corridors arising from $t$--edges not lying in the boundary
of any 2--cell). 
 In the corridor case, $e_1$ and $e_r$ are called the \emph{ends} of the corridor.  

Let $c_1, \ldots, c_{r-1}$ be the 2--cells
of $C$ numbered so that $e_i$ and $e_{i+1}$ are part of $\partial
c_i$ for all $i$.  Refer to one vertex of $e_i$ as the left and the other as the
right, depending on where it lies as we travel along $\sigma$. The
\emph{left} (\emph{right}) \emph{side} of $C$ is the edge--path in
$\Delta$ that follows the boundary--cycle $\partial c_1$ from the left (right) vertex
of $e_1$ to the left (right) vertex of $e_2$, and then $\partial
c_2$ from the left (right) vertex of $e_2$ to the left (right)
vertex of $e_3$, and so on, terminating at the left (right) vertex
of $e_r$. Note that the sides of $C$ need not be embedded paths in
$\Delta^{(1)}$.  In the case of the ring, orienting $\sigma$
anticlockwise in $\Delta$, we call the left side is the
\emph{inside} and the right side is the \emph{outside}.

The following lemma contains the basic, crucial observations about $t$--rings and $t$--corridors. 

\begin{lemma}
Suppose each $r \in \RR$ contains precisely zero or two letters
$t^{\pm 1}$.  Then every $t$--edge of $\Delta$ lies in  either a
$t$--corridor or  a $t$--ring, and the interiors of distinct $t$--rings/corridors 
are disjoint.
\end{lemma}

In the case of the presentation $\PP_k$ of Section~\ref{Phi} it will be profitable to have more general definitions:  we define an \emph{$s_k\,t$--corridor} or \emph{$s_k\,t$--ring} in a
$\PP_k$--van~Kampen diagram $\Delta$ by 
allowing the duals of edges
in $\sigma$ to be $s_k$-- or $t$--edges and allowing $\sigma$ to be a
simple edge--path whose initial and terminal vertices are either
${v_{\infty}}$ or dual to 2--cells with boundary label
$f^{-1}s_{k-1}f{s_k}^{-1}{s_{k-1}}^{-1}$. So an end
of an $s_k\,t$--corridor is either in $\partial \Delta$ or in the
boundary of a 2--cell labelled $f^{-1}s_{k-1}f{s_k}^{-1}{s_{k-1}}^{-1}$.

\section{Salient features of $\Phi_k$} \label{Phi_k salt}

\subsection{How to shortcut powers of $t$ in $\Phi_k$}
\label{shortcutting}

We fix a positive integer $k$.
The following proposition concerns the existence of \emph{shortcut diagrams} that distort $\langle t \rangle$ in
$\PP_k$.

\begin{prop} \label{shortcut diagram} There is a constant $C$, depending only on $k$, such that for all $m>0$ there is
a word $u_m$ of length $\ell(u_m) \leq C m^{1/k} + C$ and a van~Kampen diagram 
$\Delta_m$ for
${t}^m{u_m}^{-1}$ over $\PP_k$  with $\IDiam(\Delta_m) \leq C m + C$.
\end{prop}

Before proving this proposition we establish a similar result about the distortion of $\langle s_k \rangle$ in
$\OO_k$.  

\begin{lemma} \label{s_k distortion}
There is a constant $C$, depending only on $k$, such that for all $m>0$ there is a word
$v_m$ of length $\ell(v_m) \leq C m^{1/k} + C$ and
a van~Kampen diagram $\Sigma_m$ for ${s_k}^m{v_m}^{-1}$ over $\OO_k$
with every vertex of $\Sigma_m$ a distance \emph{(}in the 1--skeleton of  ${\Sigma_m}$\emph{)} at most 
$C m + C$ from the portion of $\partial \Sigma_m$ labelled  $v_m$.
\end{lemma}

\begin{figure}[ht]
\psfrag{v}{$v_m$}%
\psfrag{s_1^n}{${s_1}^n$}%
\psfrag{D_nm}{$D_{nm}$}%
\psfrag{D_n}{$D_n$}%
\psfrag{rho}{$\rho$}%
\psfrag{f^n}{$f^n$}%
\psfrag{www}{$w_n$}%
\psfrag{w_nm}{$w_{nm}$}%
\psfrag{v}{$v_m$}%
\psfrag{Sigma_m}{$\Sigma_m$}%
\psfrag{s_k^m}{${s_k}^m$}%
\psfrag{hatD_nm}{$\hat{D}_{nm}$}%
\psfrag{hatw_nm}{$\hat{w}_{nm}$}%
\centerline{\epsfig{file=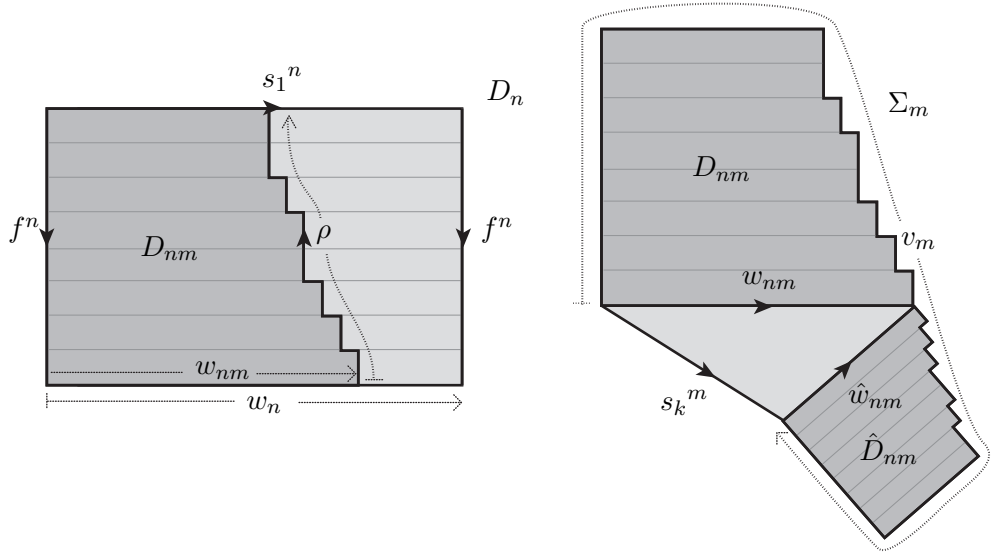}} \caption{The $\OO_k$--van~Kampen
diagrams $D_n$ and $\Sigma_m$.} \label{D_n}
\end{figure}

\Proof Lemma~\ref{binomial coefficients lemma} tells us that over
the sub--presentation
$$\langle \, s_1,
 \ldots, s_k, f \mid f^{-1} s_k f = s_k\,; \ \forall i <k, \,
f^{-1} s_i f = s_i s_{i+1} \, \rangle$$ of $\OO_k$, the word
$f^{-n} {s_1}^n f^n$ equals a (positive) word $w_n \in \set{s_1,
\ldots, s_k}^{\sstar}$ such that $w_n=\prod_{i=1}^{k} {s_i}^{n
\TB{n}{i-1}}$ in $\langle s_1, \ldots, s_k \rangle\cong \Z^k$. Fix $m>0$ and let $n$ be
the least integer such that $n\SB{n}{k-1} \geq m$. Then $(n-k+1)^k \ < \ (k-1)!\,m$ because
$$(n-1) \B{n-1}{k-1} \ = \ \frac {1}{(k-1)!} (n-1)^2(n-2) \ldots (n-k+1) \  < \  m.$$ And so 
\begin{eqnarray}
n \ < \ {(k-1)!}^{1/k}\,m^{1/k}+k-1 \ < \ k\,m^{1/k}+k. \label{n<}
\end{eqnarray}
Let $D_n$ be the
van~Kampen diagram for $f^{-n} {s_1}^n f^n \, {w_n}^{-1}$ obtained by stacking $f$--corridors $n$ high, embedded in the plane as illustrated in Figure~\ref{D_n}. Let
$w_{nm}$ be the shortest prefix of $w_n$ in which the letter
${s_k}$ occurs $m$ times.  Let $\rho$ be the edge--path from
the vertex of $\partial D_n$ at which $w_{nm}$ ends, to the portion of $\partial D_n$ labelled by $f^{-n} {s_1}^n  f^n$ , that
proceeds by travelling up where possible and left otherwise.  As $\rho$ never travels left twice consecutively, its length is at most $2n$.  
Cutting $D_n$ along $\rho$ we obtain two diagrams; one, $D_{nm}$ (which appears shaded in Figure~\ref{D_n}), shows
that $w_{nm}$ is equal in $\Theta_k$ to  a word of length at most $4n$.

Next we use the relations $[s_i,s_k]$ to
gather all the letters  $s_k$ in $w_{nm}$ to the left and produce
a word ${s_k}^m \hat{w}_{nm}$. Diagrammatically, this is
achieved by attaching $s_k$--corridors to $D_{nm}$.  (These fill the triangular region in Figure~\ref{D_n}.)  Then we use the
second stable letter $g$ of $\Theta_k$, which acts on $\langle s_1, \ldots, s_k
\rangle$ thus:
$$g^{-1} s_{k-1} g = s_{k-1}, \ \ g^{-1} s_k g = s_k \ \ \text{ and } \ \
g^{-1} s_i g = s_i s_{i+1}, \forall i <k-1\,.$$  The word $g^{-n}
{s_1}^n  g^n$ and the word obtained from $w_n$ by removing all letters $s_k$ represent the same element in $F_{k-1} \rtimes \langle g \rangle$.  This equality is exhibited by a van~Kampen diagram $\hat{D}_n$ obtained by retracting $D_n$ in the obvious manner, and a subdiagram $\hat{D}_{nm}$ (a
retraction of $D_{nm}$) of $\hat{D}_n$ portrays an equality in
$F_{k-1} \rtimes \langle g \rangle$ of $\hat{w}_{nm}$ with a word
of length at most $4n$.  Attach $\hat{D}_{nm}$ along
$\hat{w}_{nm}$ to get a diagram $\Sigma_m$ short--cutting ${s_k}^m$
to a word $v_m$ of length at most $8 n$ as shown in Figure~\ref{D_n}.
Because $D_{nm}$ and $\hat{D}_{nm}$ are stacks of $O(m^{1/k})$ corridors, it is possible to reach the portion of $\partial\Sigma_m$ labelled by $v_m$ from any vertex of $D_{nm}$ or $\hat{D}_{nm}$ by traversing  $O(m^{1/k})$ edges of $\Sigma^{(1)}$.
  The triangular region is a stack of $O(m)$ corridors which 
can be crossed to reach $\partial D_{nm}$.  
The assertion of the lemma then follows from the hypothesis 
that $\ell(v_m) \leq C m^{1/k} + C$ and the inequality (\ref{n<}).\qed

\begin{figure}[ht]
\psfrag{b}{$b$}%
\psfrag{t^m}{$t^m$}%
\psfrag{u}{$u_m$}%
\psfrag{s_k^m}{${s_k}^m$}%
\psfrag{h}{$\mbox{$\hat{s}_k$}^m$}%
\psfrag{S}{$\Sigma_m$}%
\psfrag{hS}{$\hat{\Sigma}_m$}%
\psfrag{D}{$\Delta_m$}%
\psfrag{a^3^m}{$b^{3^m}$}%
\centerline{\epsfig{file=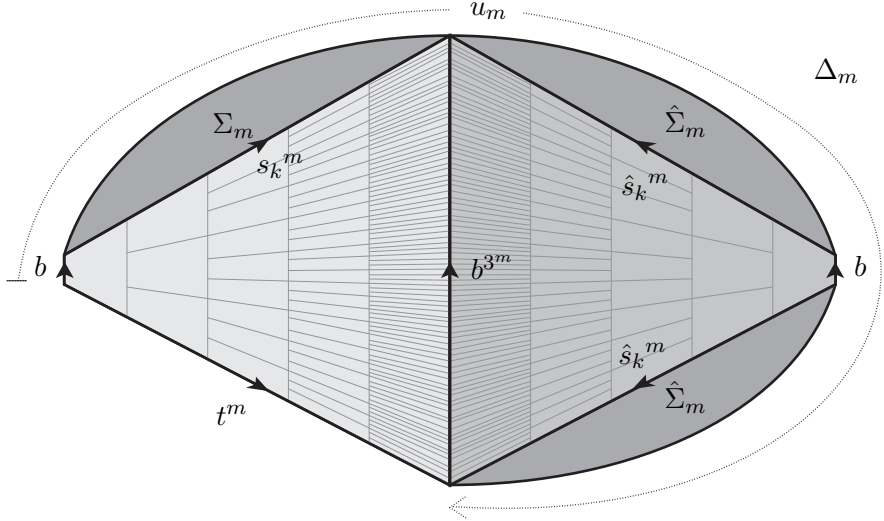}} \caption{The
$\PP_k$--van~Kampen diagram $\Delta_m$ shortcutting $t^m$.}
\label{Delta_m}
\end{figure}

\ms \ni \emph{Proof of Proposition~\ref{shortcut diagram}.}  
Using stacks of $s_k\,t$--corridors and $\hat{s}_k$--corridors in the obvious way,
one can construct van~Kampen diagrams demonstrating that 
$t^{-m} b {s_k}^m = b^{3^m}$ and $\mbox{$\hat{s}_k$}^{-m} b
\mbox{$\hat{s}_k$}^{m}=b^{3^m}$ over the subpresentations
$\langle b, s_k , t \mid t^{-1} b s_k = b^3, \ {s_k}^{-1} b t = b^3
\rangle$   and $\langle b,
\hat{s}_k \mid \mbox{$\hat{s}_k$}^{-1} b \hat{s}_k = b^3 \rangle$ of $\PP_k$, respectively.
Join these diagrams to give a diagram that demonstrates that $\mbox{$\hat{s}_k
$}^{-m}b \mbox{$\hat{s}_k$}^m = t^{-m} b {s_k}^m$ in $\Phi_k$ and that has intrinsic
diameter at most a constant times $m$. Then obtain $\Delta_m$ as
shown in Figure~\ref{Delta_m}: attach the $\OO_k$--van~Kampen
diagram $\Sigma_m$ of Lemma~\ref{s_k distortion} along ${s_k}^m$, and attach a copy of the
corresponding $\hat{\OO}_k$--van~Kampen diagram $\hat{\Sigma}_m$ and its mirror--image
along $\mbox{$\hat{s}_k$}^{m}$ and $\mbox{$\hat{s}_k$}^{-m}$. 

The asserted
bound on $\IDiam_{\PP_k}(\Delta_m)$ holds for the following reasons. Every vertex in the $t^{-m} b {s_k}^m = b^{3^m}$ and $\mbox{$\hat{s}_k$}^{-m} b
\mbox{$\hat{s}_k$}^{m}=b^{3^m}$ subdiagrams is a distance $O(m)$ from  $\Sigma_m$ or from one of the copies of $\hat{\Sigma}_m$, or from the portion of $\Delta_m$ labelled  $t^m$.  The distance from any vertex of a subdiagram $\Sigma_m$ or $\hat{\Sigma}_m$ to $\partial \Delta_m$ is $O(m)$ by Lemma~\ref{s_k distortion}.   Thus one can reach $\partial \Delta_m$ from any of its vertices  within the claimed bound.  As $\ell(\partial \Delta_m)=O(m)$, one can then follow the boundary circuit to the base point. \qed

\begin{rem} \label{shortcut area} The {\em area} of 
the diagram $\Delta_m$   in the above proof is exponential in $m$.
\end{rem} 

\subsection{Diagrams short--cutting powers $t$ in $\Phi_k$ are \emph{fat}} \label{are fat}

A key feature of the diagram  $\Delta_m$ constructed in the previous section is that it contains a path labelled by a large power $b^{3^m}$ of $b$.  In the following proposition we will show that such behaviour is common to \emph{all} $\PP_k$--diagrams $\Delta$ that greatly shortcut $t^m$ and we will deduce that all such $\Delta$ have large extrinsic (and hence intrinsic) diameter --- that is, are \emph{fat}.  
This means that such \emph{shortcut diagrams} cannot be inserted to significantly reduce the intrinsic diameter of a van~Kampen diagram.    

For a word $w$ and letter $t$, we denote the exponent sum of letters  $t$ in $w$ by $h_t(w)$.

\begin{prop} \label{shortcuts are fat} There exists a constant $C$, depending
only on $k$, with the following property: if
$u$ and $w$ are words such that $u = w$ in $\PP_k$ and $w \in \set{t,
t^{-1}}^{\sstar}$, and $\Delta$ is a 
$\PP_k$--van~Kampen diagram for $wu^{-1}$, homeomorphic to
a 2--disc, then 
$$ \IDiam(\Delta) \ \geq \  \EDiam_{\PP_{k}}(\Delta) \ \geq  \ \frac{C\abs{h_t(w)}}{1+\ell(u)}.$$ \end{prop}

\begin{figure}[ht]
\psfrag{u}{$u$}
\psfrag{u_0}{$u_0$}
\psfrag{u_1}{$u_1$}
\psfrag{wi}{$w_i$}
\psfrag{w}{$w$}
\psfrag{w_0}{$w_0$}
\psfrag{w_1}{$w_1$}
\psfrag{w_2}{$w_2$}
\psfrag{A}{$A_i$}
\psfrag{f}{$f$}
\psfrag{g}{$g$}
\psfrag{s1}{$\sstar$}
\psfrag{s2}{$\sstar'$}
\psfrag{D}{$\Delta$}
\psfrag{D_0}{$\Delta_0$}
\psfrag{D_1}{$\Delta_1$}
\psfrag{C_1}{$C_1$}
\psfrag{C_2}{$C_2$}
\centerline{\epsfig{file=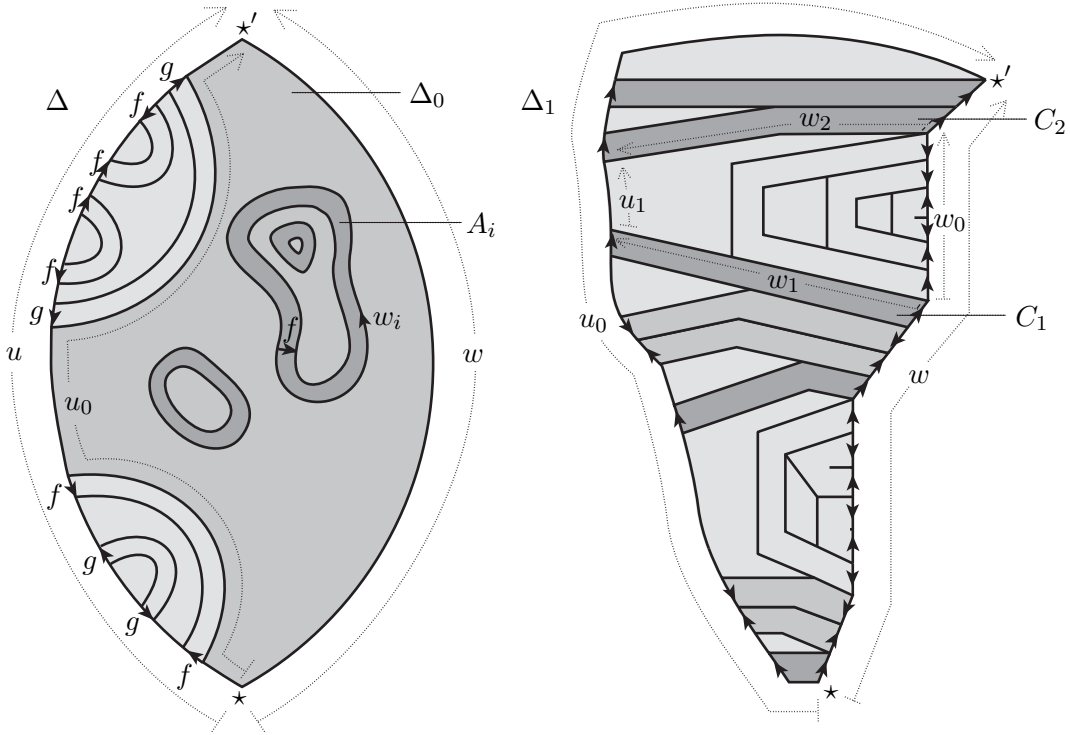}} \caption{Corridors in the
van~Kampen diagrams $\Delta$ and $\Delta_1$.}
\label{stack}
\end{figure}

\Proof  We may assume, without loss of generality, that $h_t(w) \geq 0$. Let $\gamma_w$ and $\gamma_u$
be the edge--paths in $\partial \Delta$ labelled  $w$ and
$u$, respectively,  such that the anticlockwise boundary circuit
$\partial \Delta$ is $\gamma_w$ followed by ${\gamma_u}^{-1}$. Let
$\sstar$ and $\sstar'$ be the initial
and terminal vertices of $\gamma_w$, respectively --- see the left diagram of Figure~\ref{stack}.

The essential idea in this proof is straightforward:
the $s_k\,t$--corridors (defined in Section~\ref{corridors and rings})  emanating from $\gamma_w$ stack up with the length of the corridors growing to roughly $3^{h_t(w)}$ at the top of the stack; taking the logarithm of this length
 gives a lower bound on extrinsic diameter of $\Delta$.
 However, fleshing this idea out into a rigorous proof
requires considerable care.  

The first complication we face
is that $s_k\,t$--corridors need not run right across $\Delta$, but can terminate at an edge in the side of an $f$-- or $g$-- corridor or ring.  
As every letter in $w$ is $t^{\pm1}$, the edges of $\partial\Delta$ connected in pairs
by $f$-- and $g$--corridors all lie in $\gamma_u$. Thus $\gamma_w$ lies in a
single connected component $\Delta_0$ of the planar complex obtained by
deleting from  $\Delta$   the interiors and ends of all the $f$-- and $g$--corridors.
 (In the left diagram of Figure~\ref{stack},   $\Delta_0$ is shaded.)  Note that  $\EDiam_{\PP_{k}}(\Delta) \geq \EDiam_{\PP_{k}}(\Delta_0)$. Thus
 it suffices to prove that $\EDiam_{\PP_{k}}(\Delta) \geq C\abs{h_t(w)}/(1+\ell(u)\,\!)$.
 We will do so by considering a further $\PP_k$--van~Kampen diagram $\Delta_1$ obtained from $\Delta_0$ by removing the interiors of a collection of $2$--disc subdiagrams $D_i$ and gluing in replacement diagrams $E_i$.  

Let $\gamma_{u_0}$ be the edge--path in $\partial \Delta_0$ from $\sstar$ to $\sstar'$ such that $\gamma_w$ followed by ${\gamma_{u_0}}^{-1}$ is the
anticlockwise boundary circuit of $\Delta_0$.

Suppose $A_i$ is an $f$-- or $g$--ring in $\Delta_0$ that is
not enclosed by another $f$-- or $g$--ring. The outer boundary circuit
$\gamma_i$ of $A_i$ is labelled by a word $w_i$ in $\set{ {s_1}^{\pm
1}, \ldots, {s_k}^{\pm 1}}^{\sstar}$.   
The retraction $\psi:\Phi_k \onto \Theta_k$ of Lemma~\ref{two retracts} maps $w_i$ to itself and satisfies the hypotheses of Lemma~\ref{distance decreasing}.  So $\psi$ induces a distance decreasing singular combinatorial map from $D_i$, the
diagram enclosed by $\gamma_i$,
 onto an $\OO_k$--van~Kampen diagram.  Therefore, without increasing extrinsic diameter, we may assume $D_i$ to be an $\OO_k$--diagram.   
By Lemma~\ref{abelian subpresentation}, $w_i$ is null--homotopic in the subpresentation
$\langle s_1, \ldots, s_k \mid \forall i \neq l, \ [s_i, s_l]
\rangle$ of $\PP_k$.

Let $E_i$ be a topological 2--disc van~Kampen diagram for $w_i$ in $\langle s_1, \ldots, s_k \mid \forall i \neq j, \ [s_i, s_j] \rangle$.

We will show that there is an $s_k\,t$--corridor in $\Delta_1$ with one side labelled by a word $w_T$ such that $\abs{h_b(w_T)}$ is large.  Recall that each $s_k\,t$--corridor in a $\PP_k$--van~Kampen diagram connects two edges that are either in the boundary of the
diagram or the boundary of a 2--cell with label  $f^{-1}s_{k-1}f{s_k}^{-1}{s_{k-1}}^{-1}$.  But every such 2--cell is part of an $f$--corridor, and as $\Delta_1$ contains no $f$-- or
$g$--corridor, each of its $s_k\,t$--corridors $C$ connects
two edges in $\partial \Delta_1$.
Each $t$--edge in $\gamma_w$ is part of an
$s_k\,t$--corridor that must either return to $\gamma_w$ or end on
$\gamma_{u_0}$.  An $s_k\,t$--corridor of
the former type connects two oppositely oriented $t$--edges in $\gamma_w$.  Consider $\Delta_1$ with $\gamma_{u_0}$ and $\gamma_w$ running
upwards with $\gamma_{u_0}$ on the left and $\gamma_w$ on the
right (as in Figure~\ref{stack}).  So $s_k\,t$--corridors of the latter type are
\emph{horizontal}, stacked one above another.

Define a horizontal corridor to be an \emph{up}-- or
\emph{down}--corridor according to whether the edge it meets on
$\gamma_w$ is oriented upwards or downwards.  And call an
up--corridor a  \emph{last}--up--corridor when it meets an edge on
$\gamma_w$ that is labelled by the final letter of a prefix $w_1$
of $w$ with the property that $h_t(w_2) \geq h_t(w_1)$ for all
prefixes $w_2$ of $w$ with $\ell(w_2) \geq \ell(w_1)$. Note that
there are exactly $h_t(w)$ last--up--corridors in $\Delta_1$.  In the right diagram of Figure~\ref{stack} the up-- and down--corridors are shaded and the last--up--corridors are darker.  That figure also depicts the scenario addressed in the following lemma.   

\begin{lemma} \label{corridors lemma}
Suppose that $C_1$ and $C_2$ are two horizontal corridors running
from $\gamma_w$ to $\gamma_{u_0}$ in $\Delta_1$, that $C_1$ is
below $C_2$, and that there are no horizontal corridors between
$C_1$ and $C_2$. Assume that the subarc $\gamma_{u_1}$ of
$\gamma_{u_0}$ connecting, but not including, the two edges where
$C_1$ and $C_2$ meet $\gamma_{u_0}$ is labelled by a word $u_1 \in
\set{{s_1}^{\pm 1}, \ldots, {s_{k}}^{\pm 1}}^{\sstar}$.  Let $w_1$
and $w_2$ be the words read right to left along the top of $C_1$
and along the bottom of $C_2$, respectively. Then we have equality
of the exponential sums $h_b(w_1)=
h_b(w_2)$.
\end{lemma}

\ni \emph{Proof of Lemma~\ref{corridors lemma}.} Let
$\gamma_{w_0}$  be the subarc of $\gamma_w$ connecting, but not
including, the two edges where $C_1$ and $C_2$ meet $\gamma_w$.
Let $w_0$ be the subword of $w$ that we read along $\gamma_{w_0}$.

Edges labelled by ${s_k}^{\pm 1}$ in $\gamma_{u_1}$ and by $t^{\pm
1}$ in $\gamma_{w_0}$ are the start of $s_k\,t$--corridors that must
return to $\gamma_{u_1}$ or $\gamma_{w_0}$, respectively, because
there are no horizontal corridors between $C_1$ and $C_2$. So
$h_{s_k}(u_1)=h_t(w_0)=0$ and $w_0$ freely reduces to the empty
word.  Moreover in $\PP_k$ we find that $u_1= u_2$ where $u_2$ is obtained from $u_1$ by
removing all occurrences ${s_k}^{\pm 1}$. Also,
 $w_1 u_2 {w_2}^{-1}$ is null--homotopic.

Now $w_1,w_2 \in \set{b^{\pm 1}, {s_1}^{\pm 1}, \ldots,
{s_{k-1}}^{\pm 1}}^{\sstar}$ as they are the labels of the sides of
$s_k\,t$--corridors and $u_2 \in \set{{s_1}^{\pm 1}, \ldots,
{s_{k-1}}^{\pm 1}}^{\sstar}$.  So by Lemma~\ref{null--homotopic
lemma}, $w_1 u_2 {w_2}^{-1}$ is null--homotopic in $$\langle b, s_1, \ldots, s_{k-1} \mid
\forall i \neq j, \ [s_i, s_j] \rangle \ \cong \ \Z \ast
\Z^{k-1}.$$ As $b^{\pm 1}$ does not occur in $u_2$, we deduce
that $h_b(w_1)= h_b(w_2)$. \qed

\ms 
Returning to our proof of Proposition~\ref{shortcuts are fat}, we note that edges that are part of $\gamma_{u_0}$ but not of $\gamma_u$ are
labelled by letters in $\set{{s_1}^{\pm 1}, \ldots, {s_{k}}^{\pm 1}}$
because they are in the sides of $f$-- or $g$--corridors.
So Lemma~\ref{corridors lemma} applies to adjacent horizontal corridors
that meet $\gamma_{u_0}$ at edges connected by an arc of
$\gamma_{u_0}$ that does not include edges from $\gamma_u$. By the
pigeonhole principle, there must be a stack $S$ of  $H :=
h_t(w) / (1+\ell(u))$ last--up--corridors that all meet a fixed subarc
of $\gamma_{u_0}$ that includes no edges from $\gamma_u$. (The
stack may also include horizontal corridors that are not
last--up--corridors.)

If $C$ is a horizontal $s_k\,t$--corridor and $c_1$ and $c_2$ are
the words along the top and bottom sides of $C$, then $\abs{h_b(c_1)} =
3 \abs{h_b(c_2)}$ when $C$ is an up--corridor, and $3\abs{h_b(c_1)}
= \abs{h_b(c_2)}$ when $C$ is a down--corridor.  Let $K \neq 0$ be the exponent sum of the letters $b$
in the word along the bottom edge of the lowest last--up--corridor
in $S$.  Note that $K$ is non--zero because the $s_k\,t$--corridor's two ends are labelled differently.  The number of up--corridors minus the number of down
corridors in $S$ is $H$ and so by Lemma~\ref{corridors lemma},
defining $w_T$ to be the word we read along the edge--path
$\gamma_{w_T}$ along the top of the highest corridor in $S$, we find
$h_b(w_T) = K\,3^H$.

Suppose $w'$ and $w''$ are prefixes of $w_T$.  If $w' = w''$ in $\PP_k$ then $w'  = w''$ in $\langle b, s_1, \ldots, s_{k-1} \mid \forall i \neq j, \ [s_i, s_j] \rangle$ by
Lemma~\ref{null--homotopic lemma}, and so $h_b(w') = h_b(w'')$. So
if $h_b(w') \neq h_b(w'')$, then the vertices $v'$, $v''$ of
$\gamma_{w_T}$ reached after
reading $w'$ and $w''$ map to different vertices in $\Cay^1(\PP_{k})$. So, because $h_a(w_T) =
K\,3^H$ and $\abs{K}\geq 1$, the number of vertices in
the image of $\Delta_1$ in $\Cay^1(\PP_{k})$ is at least $3^H$.  But we must say more:
there are $3^H$ prefixes $w'$ of $w_T$ that all
end $b^{\pm 1}$ and for which $h_b(w')$ are all
different.  The vertices $v'$ at the end of the arc labelled by these $w'$ are not in the interior of any
of the diagrams
 $E_i$ defined prior to Lemma~\ref{corridors lemma} 
 (as $E_i$ contain no $b$--edges) --- so the number of vertices in the
image of $\Delta_0 \ssm \left(\bigcup \textup{Int} D_i \right) =  \Delta_1 \ssm \left(\bigcup \textup{Int} E_i \right)$ in $\Cay^1(\PP_{k})$
is at least $3^H$.

\ms

The number of vertices in a closed ball $B(r)$ of radius $r$ in
$\Cay^1(\SS_{k,m})$ is at most $c^r$ for some constant $c$ depending
only on the valence of the vertices in 
$\Cay^1(\SS_{k,m})$, and hence on the number of defining
generators in $\SS_{k,m}$. Thus
$$\EDiam_{\PP_k}(\Delta) \ \geq \ \EDiam_{\PP_k}(\Delta_0) \ \geq \ \log_{c} 3^H \ = \   \frac{h_t(w)}{1+\ell(u)}\log_{c} 3.$$
\qed

\subsection{An upper bound for the distortion of $\langle t \rangle$ in $\Phi_k$} \label{d and d1}

One conclusion of Proposition~\ref{shortcut diagram} was that it is possible to express
 $t^m\in\Phi_k$ as a word of length $\preceq m^{1/k}$.  The main result of this section, Proposition~\ref{t-distortion}, is that no greater distortion of $\langle t \rangle$ is possible;
this will be crucial in Section~\ref{main thm section} when we come to  analyse $\SS_{k,m}$--van~Kampen diagrams by breaking them down into $\PP_k$-- and $\QQ_m$--subdiagrams meeting along arcs labelled by powers of $t$.

\begin{lemma}  \label{Distortion} If $u_0$ is a word in the generators of $\Theta_k$
that equals ${s_k}^m$ in the group, then $$\abs{m} \leq
\sum_{i=1}^{k} \B{c+k-1-i}{c-1} \ell_{s_i}(u_0),$$ where
$c:=(\ell_f(u_0)+\ell_g(u_0))/2$.
\end{lemma}

\Proof If we represent elements of $ \langle s_1, \ldots, s_k
\mid\, \forall \, i,j,\ [s_i,s_j]  \rangle\cong \Z^k $ as column--vectors,
then the actions of $f$ and $f^{-1}$ by conjugation are given by left multiplication by the
$k \times k$ matrices
$$\left(%
 \begin{array}{cccccc}
  1 &  &  &  &  \\
  1 & 1 &  &  &  \\
   & 1 & \ddots &  &  \\
   &  &   \ddots & 1  &  \\
   &  &   & 1 & 1 \\
\end{array}%
\right) \ \text{ and } \ \left(%
 \begin{array}{cccccc}
  1 &  &  &  &  \\
  -1 & 1 &  &  &  \\
   1  & -1 &  \ddots &  &  \\
   -1  & 1 &  \ddots & 1  &  \\
   \vdots &   \vdots &  \ddots & -1 & 1 \\
\end{array}%
\right),$$ respectively.  Similarly, we can give matrices for
the actions of $g$ and $g^{-1}$.

We will inductively obtain words $u_1, \ldots, u_c$, all of which
equal ${s_k}^m$ in $\OO_k$, as follows.  In every van~Kampen diagram
for $u_i{s_k}^{-m}$, the letters $f^{\pm 1}, g^{\pm 1}$ in $u_i$
occur in pairs $f,f^{-1}$ or $g,g^{-1}$ connected by $f$--corridors
and $g$--corridors. 
The geometry of these corridors necessitates that $u_i$ 
have a subword of the form $f^{\mp 1} v f^{\pm 1}$ or $g^{\mp
1} v g^{\pm 1}$, where $v \in \set{{s_1}^{\pm 1}, \ldots,
{s_k}^{\pm 1}}\*$.  Obtain $u_{i+1}$ from $u_i$ by replacing every
letter $s_j$ ($j =1,2, \ldots, k$) in $v$ by a word of minimal length in $\Z^k$ that
equals $f^{\mp 1} s_j f^{\pm 1}$ or $g^{\mp 1} s_j g^{\pm 1}$, as
appropriate. Let $A$ be the $k \times k$
matrix with ones in every entry on and below the diagonal and
zeros elsewhere.
 Comparing $A$ with the four matrices discussed above we see that $\ell_{s_i}(u_{i+1})$ is at most the $i$--th entry in the column vector $A \, (\ell_{s_1}(u_0), \ldots,
\ell_{s_k}(u_0) )^{\textup{tr}}$. 

Using the identity $$\sum_{n=0}^r
\B{n}{j} = \B{r+1}{j+1},$$ we calculate that $$A^c = \left(%
 \begin{array}{cccccc}
  \TB{c-1}{c-1} &  &  &  &  \\
  \TB{c}{c-1} & \TB{c-1}{c-1} &  &  &  \\
   \vdots  & \TB{c}{c-1} &  \ddots &  &  \\
   \TB{c+k-3}{c-1}  &  \vdots &  \ddots & \TB{c-1}{c-1}  &  \\
   \TB{c+k-2}{c-1} &   \TB{c+k-3}{c-1} &  \hdots & \TB{c}{c-1} & \TB{c-1}{c-1} \\
\end{array}%
\right).$$ Now $\abs{m} = \ell_{s_k}(u_c)$ and so is at most the
$k$--th entry in  $A^c \, (\ell_{s_1}(u_0), \ldots,
\ell_{s_k}(u_0)\,)^{\textup{tr}}$.  The asserted bound follows.
\qed

\ms Recall that $h_t(w)$ denotes the exponent sum of letters $t^{\pm1}$ in $w$.  

\begin{prop} \label{t-distortion}
Suppose $w$ is a word equalling $t^m$ in $\PP_k$.  Let $n:= \ell(w) - \ell_{t}(w)$. Then
$$\abs{m} \leq K\,n^k + h_t(w),$$ where $K$ is a constant
depending only on $k$.
\end{prop}

\Proof The retraction $\psi :
\Phi_k \onto \Theta_k$ of Lemma~\ref{two retracts} maps $w$, letter--by--letter, to a
word $u$  equalling ${s_k}^m$ in $\Theta_k$.    Each $t$ is mapped by
$\psi$ to $s_k$ which is central in $\Theta_k$.  So if we define
$w_0$ to be the word obtained from $w$ by removing all letters
$t^{\pm 1}$ and $u_0$ to be the word obtained by applying $\psi$, letter--by--letter, to $w_0$ then
$u_0 = {s_k}^{m-h_t(w)}$ in $\OO_k$.

By Lemma~\ref{Distortion} $$\abs{m-r} \leq \sum_{i=1}^{k}
\B{n+k-1-i}{n-1} \ell_{s_i}(u_0) \leq \B{n+k-2}{n-1}n.$$ The
asserted bound then follows because, for a
suitable constant $K$ depending only on $k$, one has
 $\B{n+k-2}{n-1} n \leq K\,n^k$ for all $n\in\mathbb N$ . \qed

\subsection{The intrinsic diameter of diagrams for $\Phi_k$}

The results in this section culminate with an upper bound on the intrinsic diameter functional of $\Phi_k$.  This will be used in Section~\ref{upper bound on EDiam section}, when we establish the upper bound on the extrinsic diameter of $\Psi_{k,m}$.   

\begin{lemma} \label{s-distortion}
If $w$ is a length $n$ word in $\PP_k$ representing an
element of the subgroup $\Z^k \cong \langle s_1, \ldots, s_k \rangle$
or $\Z^k \cong \langle \hat{s}_1, \ldots, \hat{s}_k \rangle$ of
$\Phi_k$, then there is a word $u$ in $\set{{s_1}^{\pm 1}, \ldots,
{s_k}^{\pm 1}}\*$ or in $\set{\mbox{$\hat{s}_1$}^{\pm 1}, \ldots,
\mbox{$\hat{s}_k$}^{\pm 1}}\*$, respectively, such that
$w=u$ in $\Phi_k$ and $\ell(s) \leq K\,n^k$, where $K$ depends only on $k$.
\end{lemma}

\Proof The retracts of Lemma~\ref{two retracts} mapping  $\Phi_k$ onto  $\Theta_k$ and
$\hat{\Theta}_k$ mean that it is enough to prove this lemma for $\Theta_k$ instead of $\Phi_k$.
The result for $\Theta_k$ can be proved in the same manner as Lemma~\ref{Distortion}. \qed

\bs

Define 
$$\hat{\PP}_k:=\langle \, b,t, s_1, \ldots, s_k, \hat{s}_1, \ldots, \hat{s}_k
\mid  \ {t}^{-1} b s_k = b^3, \ {s_k}^{-1} b t = b^3, \ \ \ \ \
$$
$$  \ \mbox{$\hat{s}_k$}^{-1}b\hat{s}_k = b^3; \  \ \forall i \neq j, \
[s_i,s_j]=1, [\hat{s}_i,\hat{s}_j]=1 \, \rangle, $$
a subpresentation of $\PP_k$.  Note that the 2--dimensional portions of every $\hat{\PP}_k$--van~Kampen are comprised of intersecting $s_k\,t$--, $s_i$--, and $\hat{s}_j$-- rings and corridors (where $1 \leq i < k$ and $1 \leq j \leq k$).    

\begin{lemma} \label{various lemma}
If $\Delta$ is a minimal area $\hat{\PP}_k$--van~Kampen diagram for a word $w$ then
 \begin{enumerate}
                \setlength{\itemsep}{0mm}
        \item[\textup{(}\textit{i}\textup{)}] amongst the corridors and rings in $\Delta$, no two cross twice,   
        \item[\textup{(}\textit{ii}\textup{)}] $\Delta$ contains no $s_i$-- or $\hat{s}_i$--rings \emph{(}$1\leq i <k$\emph{)},
        \item[\textup{(}\textit{iii}\textup{)}] $\Delta$ contains no $\hat{s}_k$--rings,
        \item[\textup{(}\textit{iv}\textup{)}] $\Delta$ contains no $s_k\,t$--rings, and
        \item[\textup{(}\textit{v}\textup{)}]  the length of each $s_i$-- and $\hat{s}_i$--corridor \emph{(}$1 \leq i <k$\emph{)} in $\Delta$ is less than $\ell(w) /2$.  
\end{enumerate}
 \end{lemma}

\begin{figure}[ht]
\psfrag{C_i}{$C_i$}
\psfrag{C_j}{$C_j$}
\psfrag{h_i}{$\tilde{C}_i$}
\psfrag{h_j}{$\tilde{C}_j$}
\psfrag{s_i}{$s_i$}
\psfrag{s_j}{$s_j$}
\psfrag{e}{$e$}
\psfrag{p}{$e'$}
\centerline{\epsfig{file=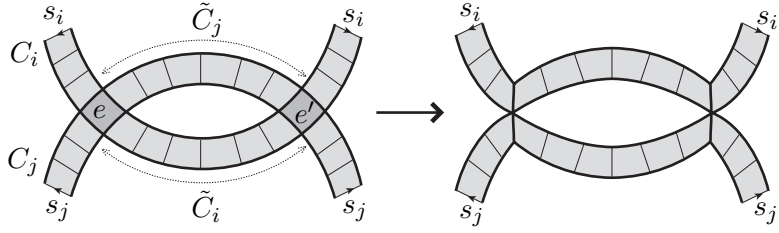}} \caption{Removing twice--crossing corridors.} \label{crossing corridors}
\end{figure}

\Proof  For (\emph{i}), first suppose that for some $i \neq j$ with $1\leq i,j<k$, there is an $s_i$--corridor or $s_i$--ring  $C_i$ that crosses an $s_j$--corridor or $s_j$--ring $C_j$ twice, intersecting at two 2--cells $e$ and $e'$.  Let $\tilde{C}_i$ and $\tilde{C}_j$ be portions of $C_i$ and $C_j$ between (but not including) $e$ and $e'$, as shown
in Figure~\ref{crossing corridors}.   By an \emph{innermost argument} we may assume that no $s_j$--corridor or $s_j$--ring intersects $\tilde{C}_i$ twice and that no  $s_i$--corridor or $s_i$--ring intersects $\tilde{C}_j$ twice.  Removing $e$ and $e'$, relabelling all the $s_i$--edges in $\tilde{C}_i$ by $s_j$ and all the $s_j$--edges in $\tilde{C}_j$ by $s_i$, and then gluing up as shown in Figure~\ref{crossing corridors}, would produce a van~Kampen diagram for $w$ of lesser area than $\Delta$.  This would contradict the minimality of the area of $\Delta$.  (Topologically, the effect of the surgery on $\Delta$ is to collapse to points arcs running through $e$ and $e'$ between opposite vertices in $\partial e$  and $\partial e'$.  
This does not spoil   planarity because no opposite pair of vertices
were identified in $\Delta$.) 

Next suppose that for some $1 \leq i <k$, an $s_k\,t$--corridor or ring $C$ crosses an $s_i$--corridor or ring $C_i$  twice at 2--cells $e$ and $e'$.  Let $\tilde{C}$ and $\tilde{C}_i$ be portions of $C$ and $C_i$ between (but not including) $e$ and $e'$.  We may assume that no $s_k\,t$--corridor or $s_k\,t$--ring intersects $\tilde{C}_i$ twice and that no $s_j$--corridor or $s_j$--ring ($1 \leq j<k$) intersects $\tilde{C}$ twice.  Let $R$ be the subdiagram between $C$ and $C_i$.  For reasons we are about to explain we may assume $\tilde{C}$ includes no 2--cells labelled by $t^{-1}bs_kb^{-3}$, (equivalently, $\partial R$ contains no $b$--edge).  This leads to  a contradiction, as above.

No 2--cell of $R$ 
 labelled $t^{-1}bs_kb^{-3}$ lies in an $s_k,t$--corridor, as such
a corridor would intersect $\tilde{C}_i$ twice.  Also $R$ includes no 2--cell that is part of an $s_k,t$--ring --- such a ring could intersect $\tilde{C}_i$ and so would have to be entirely in $R$; but then, by another \emph{innermost argument}, this $s_k,t$--ring would only contain 2--cells labelled by $t^{-1}bs_kb^{-3}$ as otherwise it would twice intersect an $s_j$--corridor for some $1\leq j <k$.  So $R$ includes no 2--cells labelled by $t^{-1}bs_kb^{-3}$, and for similar reasons, no 2--cells labelled by $\mbox{$\hat{s}_k$}^{-1}b\hat{s}_kb^{-3}$.  So all $b$--edges in $\partial R$ are identified in pairs and are amongst the 1--dimensional portions of $R$.  It follows that there are no $b$--edges in $R$, or there is a $b$--edge in $R$ that is connected to the rest of $R$ at only one vertex, or there is a 2--disc
component of $R \ssm \partial R$ that does not adjoin $C_i$.

The second case is impossible because it would imply that $\Delta$ was not reduced.
The third case cannot happen because for some $1\leq j<k$, an $s_j$--corridor would have to cross $C$, travel through this 2--disc portion and then cross $C$ again.  Thus
$R$ has no $b$--edges.

The same arguments tell us that for all $1\leq i,j \leq k$,  no $\hat{s}_i$-- and $\hat{s}_j$-- corridors or rings in $\Delta$ can cross twice.  No other combination of rings and corridors can cross even once.    

Now $(\emph{ii})$ follows immediately from $(\emph{i})$, as do $(\emph{iii})$ and $(\emph{iv})$ in the cases where the rings include 2--cells other than those labelled by $t^{-1}bs_kb^{-3}$ or $\mbox{$\hat{s}_k$}^{-1}b\hat{s}_kb^{-3}$.
In the remaining cases the outer boundary of the ring is labelled by a freely reducible word in $\set{b^{\pm 1}}^{\sstar}$, which contradicts the minimality of the area of $\Delta$.

For $(\emph{v})$, suppose $C$ is an $s_i$-- or $\hat{s}_i$--corridor.  The boundary circuit of $\Delta$ is comprised of the two ends of $C$ and two edge--paths, one of which, call it $\alpha$, must have length less than $\ell(w)/2$.  By (\emph{i})--(\emph{iv}) a different corridor connects each on a side of $C$ to $\alpha$.  So $C$ has length less than $\ell(w)/2$.      \qed

\begin{lemma}\label{subpresentation} There exists a constant $K$ such that
$\IDiam_{\hat{\PP}_k}(n)\le Kn$ for all $n\in\mathbb N$. 
\end{lemma}

\Proof Suppose $\Delta$ is a minimal area $\hat{\PP}_k$--van~Kampen diagram for a word $w$.  All $s_k\, t$--corridors in $\Delta$ are embedded: if a portion (of non--zero length) of the path along the side of an $s_k\,t$--corridor formed an edge--loop then that would have to enclose a zero area subdiagram (by the results of Lemma~\ref{various lemma}) and therefore be labelled by a non--reduced path --- but then $\Delta$ would not be a minimal area diagram as there would be an inverse pair of 2--cells on the corridor.  So from any given vertex in $\Delta$, travelling across $s_k\, t$--corridors at most $n/2$ times, we  meet either an $s_i$--corridor ($1 \leq i <k$), or an $\hat{s}_i$--corridor ($1 \leq i <k$), or $\partial \Delta$.  In the former two  cases one can follow a path of length at most $n/4$ along a side of the $s_i$-- or $\hat{s}_i$--corridor to $\partial \Delta$.  From any point on $\partial \Delta$ one can reach the base vertex by following the boundary circuit.   
\qed

\bs
 
\begin{prop} \label{intrinsic diam of PP_k} $\IDiam_{\PP_k}(n) = O(n^k)$.
\end{prop}

\Proof  Suppose $w$ is a null--homotopic word in $\PP_k$ and $\Delta$ is a van~Kampen diagram for $w$. We will construct a new
van~Kampen diagram $\Delta'$ for $w$ that satisfies the claimed
bound on intrinsic diameter. We begin the construction of $\Delta'$
by taking an  edge--circuit $\xi$ of length $\ell(w)$ in
the plane to serve as $\partial\Delta'$.  We direct and label the edges of $\xi$ so that one reads $w$ around it.

The occurrences of  $f$ and $f^{-1}$  in $w$ are paired so that the
corresponding edges of $\partial\Delta$ are joined by corridors
in $\Delta$. For each such pair, there must be a subword
$f^{-1}w_0f$ in some cyclic conjugate of $w$ or $w^{-1}$ such that both
$w_0$ and $f^{-1}w_0f$ represent elements of the subgroup $\Z^k =
\langle s_1, \ldots, s_k \rangle$.   Join the pairs of $f$-- and
$f^{-1}$--edges in $w$ by $f$--corridors running through the interior of
$\xi$. It follows from
Lemma~\ref{s-distortion} that the lengths of both sides of
these $f$--corridors are at most a constant times
$(2+\ell(w_0) )^k$.

Likewise, insert corridors into the
interior of $\xi$ joining pairs of $g$-- and $g^{-1}$--edges, of $\hat{f}$-- and $\hat{f}^{-1}$--edges, and of $\hat{g}$ and
$\hat{g}^{-1}$--edges.  There is no obstruction to planarity in the
2--complex because we are mimicking the layout
of corridors in $\Delta$.

To complete the construction of $\Delta'$, we shall fill the 2--disc holes
inside $\xi$.  The boundary circuit $\eta$ of each hole is made up
of the sides of corridors $C_1, \ldots, C_r$ and a number of disjoint subarcs of
$\xi$. The length of $\eta$ is at most $n^k$, up to a
multiplicative constant, because these disjoint subarcs contribute at
most $n$ and the lengths of the sides of $C_1, \ldots, C_r$ each contribute at
most a constant times ${n_i}^k$, where $\sum_{i=1}^r n_i \leq
n$. By Lemma~\ref{subpresentation}, these circuits can be filled
by  $\hat{\PP}_k$--van~Kampen diagrams with intrinsic diameter at most a constant
times $n^k$. And, as the corridors have length at most a constant
times $n^k$, we deduce that $\IDiam_{\PP_k}(\Delta') = O(n^k)$, as
required. \qed

\begin{rem} \label{PPk Dehn} Let us consider why
the Dehn function of $\PP_k$ is at most $n \mapsto K^{n^k}$ for some constant $K>0$.  

 Lemma~\ref{various lemma} implies that the Dehn function of $\hat{\PP}_k$ is at most ${K_1}^n$ for some $K_1 >0$: the total contribution of the 2--cells labelled by $[s_i, s_j]$ or $[\hat{s}_i, \hat{s}_j]$ for some $i \neq j$ is at most $n^2$; removing all $s_i$-- and $\hat{s}_i$--corridors \textup{(}for all $1\leq i <k$\textup{)} leaves components with linear length boundary circuits filled by minimal area van~Kampen diagrams over    
$$\langle b,t, s_k, \hat{s}_k
\mid  \ {t}^{-1} b s_k = b^3, \ {s_k}^{-1} b t = b^3, \ \mbox{$\hat{s}_k$}^{-1}b\hat{s}_k = b^3 \rangle,$$
and a standard corridors argument shows the Dehn function of this subpresentation is bounded above by an exponential function.
The construction of diagrams in Proposition~\ref{intrinsic diam of PP_k} then establishes that the  Dehn function of $\PP_k$ is at most  $n \mapsto K^{n^k}$ for some $K>0$.
\end{rem}

\section{Diameter in $\Gamma_m$} \label{d and d2}

In this section we establish an upper bound on the intrinsic, and hence extrinsic, diameter of null--homotopic words $w$ in the presentation $\QQ_m$ for $\Gamma_m$.  Also we show how the \emph{shortcuts} of Proposition~\ref{shortcut diagram} lead to an improved bound on extrinsic diameter when we regard $w$ as a word in the presentation $\SS_{k,m}$ for $\Psi_{k,m}$.


\begin{prop} \label{extrinsic diam
bound} Fix integers $n>0$ and $1< k <m$.  Suppose $w$ is a
null--homotopic word in the presentation $\QQ_m$ for $\Gamma_m$,
that $\ell(w) - \ell_t(w) \leq n$, and that $\ell_t(w) = O(n^k)$.
Then there is a $\QQ_m$--van~Kampen diagram $\Delta$ for $w$ with
$\IDiam_{\QQ_m}(\Delta) = O(n^{m+1})$. Moreover, as an
$\SS_{k,m}$--van~Kampen diagram, $$\EDiam_{\SS_{k,m}}(\Delta)=O\left(n^{\max\set{{1+\frac{m}{k}},\,k}}\right).$$
\end{prop}

To prove this result we will need some purchase on the geometry of $\tau$-- and $\sigma$--corridors.  This is provided by the following two lemmas.

\begin{lemma}\label{tau-corridors}
Suppose $\Delta$ is a minimal area $\QQ_m$--van~Kampen diagram $\Delta$ for a null--homotopic word $w$ and that $\CC$ is a $\tau$--corridor in $\Delta$.  Let $w_0$ be the word along the sides of $\CC$.  If $w_1$ is a subword of $w_0$ and $w_1 \neq t^{\pm1}$ then $w_1$ cannot represent the same group element as a non--zero power of $t$.  In particular, $\CC$ is embedded; that is, the sides of $\CC$ are simple paths in $\Delta$. 
\end{lemma}

\Proof As $\Delta$ is of minimal area, $w_0$ is freely reduced as a word in $$\set{{a_1}^{\pm 1}, \ldots, {a_{m-1}}^{\pm 1}, (a_mt)^{\pm 1}, T^{\pm 1}}^{\sstar}.$$  Killing $t$ and $\tau$ retracts $\QQ_m$ onto $$\UU_m:=\langle a_1, \ldots, a_m, \s, T \mid  \s^{-1}a_m \s = a_m, \  \forall i < m, \ \s^{-1} a_i \s = a_ia_{i+1} \rangle,$$ in which $a_1,\ldots, a_m,T$ generate a free subgroup.  This retraction sends $w_1$ to a word in $\set{{a_1}^{\pm 1}, \ldots, {a_{m-1}}^{\pm 1}, {a_m}^{\pm 1}, T^{\pm 1}}^{\sstar}$ that is freely reduced and therefore non--trivial.  The result follows. \qed

\bs 
The analogue of this result for $\sigma$--corridors is more complex.  

\begin{lemma}\label{sigma-corridors}
Suppose $\Delta$ is a minimal area $\QQ_m$--van~Kampen diagram $\Delta$ for a null--homotopic word $w$ and that $\CC$ is a $\sigma$--corridor in $\Delta$.  Then $\CC$ is a topological 2--disc subdiagram of $\Delta$ with boundary label $\sigma^{-1} w_0 \sigma {w_1}^{-1}$, where $w_0$ and $w_1$ are  freely reduced
words in $\set{{a_1}^{\pm 1}, \ldots, {a_m}^{\pm 1}}^{\sstar}$.  Moreover, defining $\rho_0$ and $\rho_1$ to be the arcs of $\partial \CC$ along which one reads $w_0$ and $w_1$, there exists $K>0$, depending only on $m$, such that from any point on $\rho_i$ \emph{(}$i=0,1$\emph{)} there is a path in $\CC^{(1)}$ of length at most $K$ to $\rho_{\abs{1-i}}$.  Also, no subword of $w_0$ or $w_1$ represents the same group element as a non--zero power of $t$.
\end{lemma}

\Proof
That $w_0$ is reduced, that $\rho_0$ is a simple path in $\Delta^{(1)}$, and that no subword of $w_0$ or $w_1$ represents the same group element as a non--zero power of $t$, are all proved as in Lemma~\ref{tau-corridors}.  And $\CC$ is a topological 2--disc by a similar argument involving the retraction $\UU_m$.

So $w_1$ is obtained from $w_0$ by replacing each ${a_i}^{\pm 1}$ by $(a_ia_{i+1})^{\pm 1}$ for all $i<m$, and then freely reducing.  The constant $K$ exists by the \emph{Bounded Cancellation Lemma} of \cite{Cooper}. 
\qed

\bs

\ni \emph{Proof of Proposition~\ref{extrinsic diam bound}.}
Let $\Delta$ be a $\QQ_m$--van~Kampen diagram for $w$ that is of minimal area.
Note that $\Delta$ contains no $t$--, $\tau$-- or $\sigma$--rings.  Removing the $\tau$-- and $\sigma$--corridors from $\Delta$ leaves a disjoint union of subdiagrams $\Delta_i$ over $$\langle a_1, \ldots, a_m, t, T \mid \forall j, [t, a_j]=1\,; \  [t,T]=1 \rangle.$$  Since there are no $t$--rings in $\Delta_i$ and the words along the sides of the $t$--corridors are reduced, from any vertex in $\Delta_i$ one can reach $\partial \Delta_i$ by following at most $\ell(\partial \Delta_i)/2$ $t$--edges across successive $t$--corridors.    

We claim that if 
$\CC$ is a $\tau$-- or $\sigma$--corridor in $\Delta$ and $u$ is the word along a side of $\CC$ then $\ell(u) \leq O(\min\set{{n_1}^m, {n_2}^m})$, where $n_1$ and $n_2$ are the lengths of the two arcs that together comprise $\partial \Delta$ and have the same end points as one side of $\CC$.   This is because killing $t$ and $\tau$ retracts $\QQ_m$ onto $\UU_m $, and the calculation used in the proof of Proposition~\ref{Distortion} 
shows that the conjugation action of the stable letter $\s$ is an  automorphism of polynomial growth of degree $m-1$.   And on  page~452 of \cite{Bridson2} the first author shows that the growth of the inverse 
automorphism is also polynomial of degree $m-1$.  

It follows that $\ell(\partial \Delta_i) \leq O(n^m)$ because the boundary of $\Delta_i$ consists of portions of $\partial \Delta$, and the sides of $\sigma$-- and $\tau$--corridors.  

We are now ready to estimate $\IDiam_{\QQ_m}(\Delta)$.  Suppose $v$ is a vertex of $\Delta$.  In the light of  Lemmas~\ref{tau-corridors} and \ref{sigma-corridors},
there is no loss of generality in assuming
$v$ is not in the interior of a $\tau$-- or $\sigma$--corridor.   Move from $v$ to $\partial \Delta$ by successively crossing $\s$-- and $\tau$--corridors as follows.  When located on a $\sigma$-- or $\tau$--corridor that has not just been crossed, follow at most $K$ (the constant of Lemma~\ref{sigma-corridors}) edges to cross to the other side; otherwise follow a maximal length embedded path of $t$--edges across some $\Delta_i$.  Finally follow $\partial \Delta$ to the base vertex $\sstar$.  Let $\rho$ denote the resulting edge--path from $v$ to $\sstar$.

It must be verified that we can indeed
 reach  $\partial \Delta$ by moving in the manner described  and that $\ell(\rho)=O(n^m)$.  First note that there are no embedded edge--loops in $\Delta$ labelled by words in $\set{t^{\pm 1}}^{\sstar}$ because the interior of such a loop could be removed and the hole glued up (as $\QQ_m$ retracts onto $\langle t \rangle$), reducing the area of $\Delta$.  Next observe that $\rho$ does not cross the same $\tau$-- or $\sigma$--corridor twice, for otherwise there would have to be an \emph{innermost} $\tau$-- or $\s$--corridor that $\rho$ crosses twice, contradicting either Lemma~\ref{tau-corridors} or Lemma~\ref{sigma-corridors}.  


So crossing the $\s$-- and $\tau$--corridors, of which there are at most $n/2$, contributes at most $Kn/2$ to the length of $\rho$.  We have already argued that each section of $\rho$ between an adjacent pair of $\s$-- or $\tau$--corridors has length $O(n^m)$.  So these sections together contribute at most $(n/2) O(n^m) = O(n^{m+1})$ to the total length.  The final section of $\rho$ is part of $\partial \Delta$ and so has length at most $O(n^{k}) < O(n^m)$.  The total is $O(n^{m+1})$ as required.   

For the bound on $\EDiam_{\SS_{k,m}}(\Delta)$, we note that  in the word metric (i.e.\ measured in the
Cayley graph)  the distance from the initial
to the terminal vertex of each arc of $\rho$ whose edges are all labelled by $t^{\pm1}$, is $O(n^{\frac{m}{k}})$ by Proposition~\ref{shortcut diagram}.  Therefore 
$\EDiam_{\SS_{k,m}}(\Delta) = O(n)O(n^{\frac{m}{k}}) = O(n^{1+\frac{m}{k}})$, as required.   \qed

\begin{rem} \label{Gamma Dehn}
The proof of Proposition~\ref{extrinsic diam bound} also shows that the Dehn function of $\Gamma_m$ is $O(n^{2m+1})$  \emph{(}cf.\ Theorem 3.4 of \textup{\cite{Bridson}}\emph{)}. 
\end{rem}

\section{Proof of the main theorem} \label{main thm section}

In this section we prove Theorem~\ref{Main thm with details} (modulo some technicalities postponed to Section~\ref{Amalgams and retractions}).    Theorem~\ref{main thm} then follows because  $k$ and $m$ can be chosen so that the ratio of the intrinsic
and extrinsic diameter filling functionals of $\Psi_{k,m}$ grows
faster than any prescribed polynomial.

Before stating the theorem we recall the definition of an \emph{alternating product expression} in an amalgam and a well known lemma --- see Lemma~6.4 in Section~III.$\Gamma$ of \cite{BrH} or Section~5.2 of \cite{Trees} .

\begin{defn} \label{alt prod decomp}
Let $\Gamma = A \ast_C B$ be an amalgam of groups $A$ and $B$
along a common subgroup $C$.  Suppose $\AA$ and $\BB$ are
generating sets for $A$ and $B$, respectively. An
\emph{alternating product expression} for $w \in (\AA^{\pm
1} \cup \BB^{\pm 1})^{\sstar}$ is a cyclic conjugate $u_1 v_1
\ldots u_p v_p$ of $w$ in which for all $i$ we have $u_i \in (\AA^{\pm
1})^{\sstar}$, $v_i \in (\BB^{\pm 1})^{\sstar}$, and if $p \neq 1$ then  neither
$u_i$ nor $v_i$ is the empty word.
\end{defn}

\begin{lemma} \label{alt prod lemma}
In the notation of Definition~\ref{alt prod decomp}, if $w \in (\AA^{\pm 1} \cup \BB^{\pm
1})^{\sstar}$ represents 1 in $\Gamma$ then in any alternating
product expression $u_1 v_1 \ldots u_p v_p$ for $w$, some $u_i$ or
$v_i$ represents an element of $C$.
\end{lemma}

\begin{thm} \label{Main thm with details}
For all integers $m >k>1$, the extrinsic and intrinsic diameter filling
functions of the group $\Psi_{k,m}$ presented by $\SS_{k,m}$, satisfy %
\begin{eqnarray*}
\EDiam_{\SS_{k,m}}(n) & = & O\!\left(n^{\max\set{1+\frac{m}{k},\,k}}\right), \\ %
n^{m/3} & = & O(\IDiam_{\SS_{k,m}}(n)).
\end{eqnarray*}
\end{thm}

\subsection{Proof of the upper bound on $\EDiam_{\SS_{k,m}}(n)$} \label{upper bound on EDiam section}

Suppose $w$ is a null--homotopic word in $\SS_{k,m}$.
Let $u_1 v_1 \ldots u_p v_p$ be an \emph{alternating product
decomposition} for $w$ where each $u_i$ and $v_i$ is a word on the generators of $\PP_k$ and
$\QQ_m$ respectively.  Take a planar edge--circuit $\eta$ around
which, after directing and labelling the edges, one reads $w$.
Decompose $\eta$ into arcs along which we read the $u_i$ and
$v_i$, and call the vertices at which these meet
\emph{alternation vertices}.  Say that an edge--path is a
\emph{$t$--arc} when it is made up of $t$--edges orientated the same way. 

Repeated appeals to Lemma~\ref{alt prod lemma}, the first of which tells us that some
$u_i$ or $v_i$ represents a word in $\set{t^{\pm 1}}\*$, allow us
to deduce that $t$--arcs  with the following properties can be inserted into the interior of
$\eta$: each $t$--arc connects two distinct alternation vertices of
$\eta$\,; any two $t$--arcs are disjoint;
and the $t$--arcs partition the interior of $\eta$ into topological
2--disc regions, whose boundary loops are labelled by words that
are null--homotopic either in $\PP_k$ or in $\QQ_m$. Call these
bounding loops
\emph{$\PP_k$--loops} and \emph{$\QQ_m$--loops}, respectively.

\begin{lemma} \label{t-arc length}
For each $t$--arc $\alpha$, let $n_{\alpha}$ be the length of the
shorter of the two subarcs of $\eta$ that share their two end
vertices with $\alpha$.  Then the length $\ell(\alpha)$ of
$\alpha$ satisfies
$$\ell(\alpha) \ \leq \ K \, {n_{\alpha}}^k$$ where $K \geq 1$ is the constant of
Proposition~\ref{t-distortion}, which depends only on $k$.
\end{lemma}

\begin{figure}[ht]
\psfrag{alpha}{$\alpha$}%
\psfrag{alpha_1}{$\alpha_1$}%
\psfrag{alpha_2}{$\alpha_2$}%
\psfrag{u_1}{$u_1$}%
\psfrag{v_1}{$v_1$}%
\psfrag{u_2}{$u_2$}%
\psfrag{v_2}{$v_2$}%
\psfrag{u_3}{$u_3$}%
\psfrag{v_3}{$v_3$}%
\psfrag{u_4}{$u_4$}%
\psfrag{v_4}{$v_4$}%
\psfrag{u_5}{$u_5$}%
\psfrag{v_5}{$v_5$}%
\psfrag{n}{$n_{\alpha}$}%
\psfrag{n_{alpha_1}}{$n_{\alpha_1}$}%
\psfrag{n_{alpha_2}}{$n_{\alpha_2}$}%
\psfrag{RR}{$\RR$}%
\psfrag{star}{$\sstar$}%
\centerline{\epsfig{file=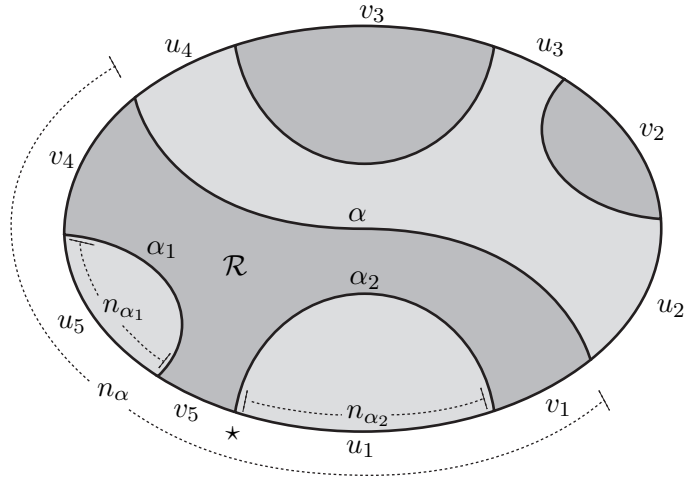}} \caption{Partitioning $\eta$
with $t$-arcs.} \label{eta}
\end{figure}

\ni \emph{Proof of Lemma~\ref{t-arc length}.} We induct on
$n_{\alpha}$. The base case, $n_{\alpha}=0$ holds vacuously.
For the inductive step, of the two regions
adjoining $\alpha$, choose $\RR$ to be that which is in the interior of
the disc bounded by $\alpha$ and the length $n_{\alpha}$ subarc of
$\eta$ --- see Figure~\ref{eta}. The boundary of $\RR$ is comprised of $\alpha$, disjoint
subarcs of $\eta$, and some $t$--arcs that we call $\alpha_1, \ldots,
\alpha_r$. Define $n_0$ to be the total length of these subarcs of
$\eta$. Then 
\begin{eqnarray}
n_{\alpha} \ =  \ n_0 + \sum_{i=1}^r n_{\alpha_i} \label{some sum
equation}
\end{eqnarray}
By induction hypothesis, $\ell(\alpha_i) \leq K \,
{n_{\alpha_i}}^k$ for all $i$. We address two cases.

\ms \ni \emph{Case: $\partial \RR$ is a $\QQ_m$--loop.} The exponent
sum of the occurrences of $t$ in any null--homotopic word in
$\QQ_m$ is zero because $\QQ_m$ retracts onto $\langle t \rangle
\cong \Z$ via the map that kills all generators other than $t$. So
$$\ell(\alpha) \ \leq \ n_0 + \sum_{i=1}^r K \, {n_{\alpha_i}}^k \ \leq \ K \,
{n_{\alpha}}^k,$$ with the second inequality following from
(\ref{some sum equation}).

\ms \ni \emph{Case: $\partial \RR$ is a $\PP_k$--loop.}
Proposition~\ref{t-distortion} and (\ref{some sum equation}) give
$$\ell(\alpha) \ \leq \ K \, {n_0}^k + \sum_{i=1}^r K \,
{n_{\alpha_i}}^k \ \leq \ K \, {n_{\alpha}}^k\;.$$

\ms \ni This completes the proof of the lemma. \qed

\ms Next we fill the $\PP_k$-- and $\QQ_m$--loops inside $\eta$ to
produce a van~Kampen diagram that satisfies the asserted bound on
extrinsic diameter. Let $n:=\ell(w)$.  By Lemma~\ref{t-arc length}, the length
of each of their boundary circuits is at most $O( n^k)$, with the contributions from portions not on $\eta$ coming entirely from   $t$--arcs.  First we fill the $\QQ_m$--loops as per Proposition~\ref{extrinsic diam bound}, 
with diagrams each of which meets $\eta$ and has extrinsic diameter $O\left(n^{\max\set{1+\frac{m}{k},\,k}}\right)$ as an $\SS_{k,m}$--van~Kampen diagram.

Next we glue a \emph{shortcut diagram} along each of the $t$--arcs as per Proposition~\ref{shortcut diagram}.  These diagrams have intrinsic diameter $O(n^k)$, measured from base vertices on $\eta$. Finally, fill the remaining $\PP_k$--loops --- all have length $O(n)$ and meet $\eta$, and by Proposition~\ref{intrinsic diam of PP_k} they can be filled by van~Kampen diagrams of intrinsic, and hence extrinsic, diameter $O(n^k)$.  

The result is an $\SS_{k,m}$--van~Kampen diagram that admits the asserted bound on extrinsic diameter. \qed

\begin{rem} \label{main Dehn fn rem}
This proof, together with Remarks~\ref{shortcut area}, \ref{PPk Dehn} and \ref{Gamma Dehn} establish that the Dehn function of $\Psi_{k,m}$ is at most $n \mapsto C^{n^k}$ for some constant $C$.   
\end{rem}

\subsection{Proof of the lower bound on $\IDiam_{\SS_{k,m}}(n)$.} \label{proof of lower bound}

\begin{figure}[ht]
\psfrag{s^n}{$\s^n$}
\psfrag{a_1^n}{${a_1}^n$}
\psfrag{tau}{$\tau$}
\psfrag{T}{$T$}
\psfrag{u}{$u$}
\psfrag{v}{$v$}
\psfrag{u_t}{$u_t$}
\psfrag{v_t}{$v_t$}
\psfrag{t^q}{$t^q$}
\psfrag{star}{$\sstar$}
\centerline{\epsfig{file=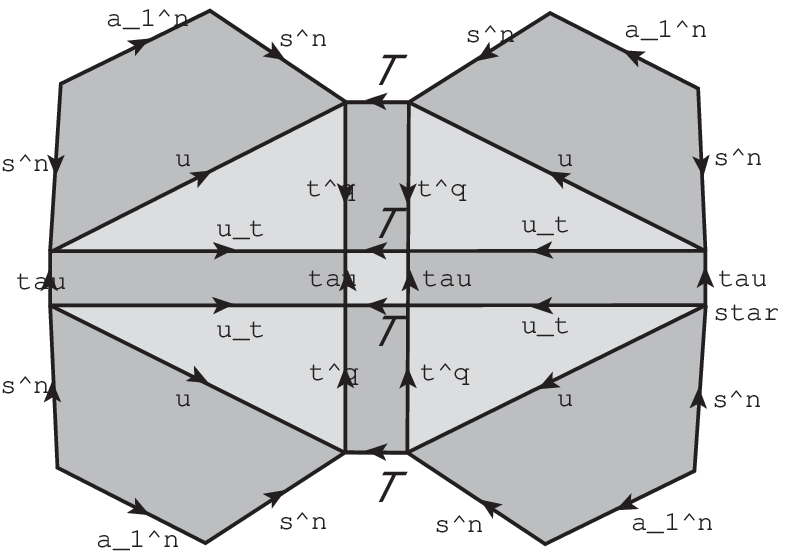}} \caption{A $\QQ_m$--van~Kampen
for $w_n$.} \label{w_n}
\end{figure}

Define $$w_n  \ := \ [\tau,
(\s^{-n}{a_1}^n\s^n)\,T\,(\s^{-n}{a_1}^{-n}\s^n)].$$
Let $u=u(a_1, \ldots, a_m)$ be the (positive) word such that 
$\s^{-n}{a_1}^n\s^n = u$ in $\QQ_m$.  Define $u_t := u(a_1, \ldots, a_{m-1}, a_mt)$.
For some $q \in \Z$ we have $ut^q = u_t$.  This plays a key role in the construction of the $\QQ_m$--van~Kampen for $w_n$, an outline for which is shown in Figure~\ref{w_n}. 

\begin{lemma}\label{extrinsic diameter lemma}
The word $w_n$ has extrinsic diameter at least $C\,n^m
-C$ in $\QQ_m$, where $C$ is a constant depending only on $m$.
\end{lemma}

\ni \emph{Proof of Lemma~\ref{extrinsic diameter lemma}.} Our
approach builds on the proof by the first author of
Theorem~3.4 in \cite{Bridson}. Suppose $\pi: \Delta \to
\Cay^2(\QQ_m)$ is a $\QQ_m$--van~Kampen diagram for $w_n$. First we
find an edge--path $\rho$ in $\Delta$, along which one reads a word
in which the exponent sum of the letters $t$ is at least $C\,n^m
-C$ for some constant $C= C(m)$. A $T$--corridor connects the two
letters $T$ in $w_n$, and along each side of this corridor we read a
word $u$ in $\set{t^{\pm 1}, \tau^{\pm 1}}^{\sstar}$ that equals
$(\s^{-n}{a_1}^{-n}\s^n)\tau^{-1}(\s^{-n}{a_1}^n\s^n)$ in $\QQ_m$.  A
$\tau$--corridor joins the $\tau^{-1}$ in this latter word to some
$\tau^{-1}$ in $u$.  Let $u_0$ be the prefix of $u$ such that the
letter immediately following $u_0$ is this $\tau^{-1}$, and then
let $\rho$ be the edge--path along the side of the $\tau$--corridor
running from the vertex at the end of $\s^{-n}{a_1}^{-n}\s^n$ to the
vertex at the end of $u_0$.  Let $v=v({a_1}, \ldots, {a_{m-1}},
{a_mt})$ be the word one reads along $\rho$. Then $u_0=
(\s^{-n}{a_1}^{-n}\s^n)v$ in $\QQ_m$. Killing $T$, $t$ and $\tau$,
retracts $\QQ_m$ onto
$$\langle a_1, \ldots, a_m, s \mid \forall \, i \neq j, \  [a_i,a_j]=1 \ ; \ \s^{-1}a_m \s = a_m, \  \forall i <
m, \ \s^{-1} a_i \s = a_ia_{i+1} \rangle,$$ in which $\bar{v} = \s^{-n}{a_1}^{n}\s^n$, where
$\bar{v}:=v({a_1}, \ldots, {a_{m-1}}, {a_m})$. By
Lemma~\ref{binomial coefficients lemma}, the exponent sum of $a_m$
in $\bar{v}$, and hence of $t$ in $v$, is $n\TB{n}{m-1}$, which is
at least $C\,n^m -C$ for some constant $C$ depending only on $m$.

Killing all generators other than $t$ defines a retraction $\phi$
of $\Gamma_m$ onto $\langle t \rangle \cong \Z$. The 0--skeleton of
$\Cay^2(\QQ_m)$ is  $\Gamma_m$ and the existence of $\rho$ shows that
image of $\phi \circ \pi: \Delta \to \Z$ has diameter at least
$C\,n^m$, since the retraction $\phi$ from
the 0--skeleton of $\Cay^2(\QQ_m)$ to $\Z$ decreases distance.
This completes the proof of the lemma. \qed

\ms Let $\Delta$ be a minimal intrinsic diameter
$\SS_{k,m}$--van~Kampen diagram for $w_n$. On account of the
retraction $\phi$, we may assume $\Delta$ to have the properties
described in Proposition~\ref{islands} (\emph{2}). That is, the
$\PP_k$--cells within $\Delta$ comprise a subcomplex whose
connected components are all simply connected unions of
topological disc subcomplexes any two of which meet at no more
than one vertex. Refer to these topological disc subcomplexes as
\emph{$\PP_k$--islands}. Around the boundary of each $\PP_k$--island
we read a word in $\set{t^{\pm 1}}\*$ that freely reduces to the
empty word because $t$ has infinite order in  $\Phi_{k,m}$.

\begin{figure}[ht]
\psfrag{t}{$t$}
\psfrag{u}{$u$}
\psfrag{v}{$v$}
\psfrag{uh}{$\hat{u}$}
\psfrag{vh}{$\hat{v}$}
\psfrag{tu}{$\theta(u)$}
\psfrag{tv}{$\theta(v)$}
\psfrag{tuh}{$\theta(\hat{u})$}
\psfrag{tvh}{$\theta(\hat{v})$}
\psfrag{I}{$I$}
\psfrag{s}{$\sstar$}
\psfrag{sb}{$\theta(\sstar)$}
\psfrag{t}{$\theta$}
\psfrag{d}{$\delta_{u,v}$}
\psfrag{hd}{$\hat{\delta}_{u,v}$}
\psfrag{gamma}{$\gamma_{u,v}$}
\psfrag{hg}{$\hat{\gamma}_{u,v}$}
\psfrag{theta}{$\theta(\partial I)$}
\psfrag{D}{$\Delta$}
\psfrag{bD}{$\bar{\Delta}$}
\centerline{\epsfig{file=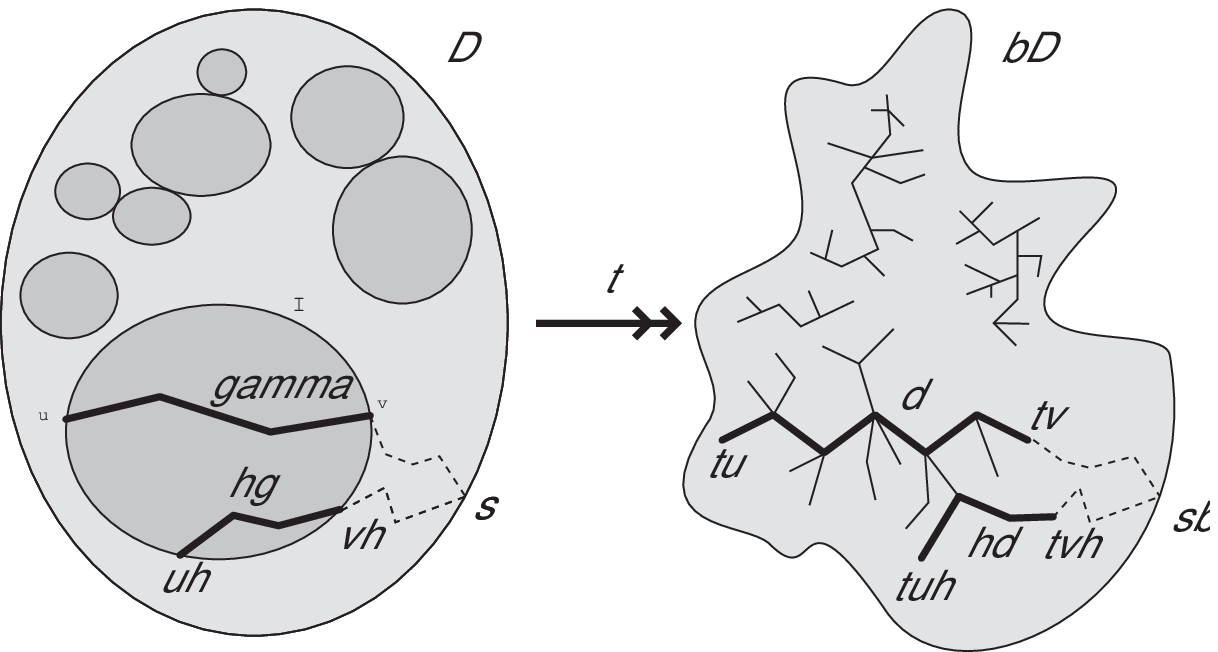}} \caption{Collapsing the $\PP_k$--islands.}
\label{collapse}
\end{figure}

\ms We obtain a $\QQ_m$--van Kampen diagram $\bar{\Delta}$ for
$w_n$ from $\Delta$ by cutting out the $\PP_k$--islands and then
gluing up the attaching cycles by identifying adjacent,
oppositely--oriented edges (i.e.\ successively cancelling pairs
$tt^{-1}$ or $t^{-1}t$ in the attaching words).  The removal of
the $\PP_k$--islands and subsequent gluing is described by a collapsing map
$\theta:\Delta \onto \bar{\Delta}$ that is injective and
combinatorial except on the $\PP_k$--islands.  This is depicted in Figure~\ref{collapse}.  Note that, no matter
what choice of $\theta$ we make, $\bar{\Delta}$ will be planar by
Lemma~\ref{collapse cut glue}, since after cutting out and gluing
up a number of the $\PP_k$--islands, the boundary circuit of every
remaining $\PP_k$--island is a simple loop.

Let $T$ be a maximal geodesic tree in the 1--skeleton of $\Delta$,
based at $\sstar$. Suppose $\gamma$ is a geodesic in $T$ from a
vertex $v$ not in the interior of a $\PP_k$--island to $\sstar$.
Define an edge--path $\bar{\gamma}$ in the 1--skeleton of
$\bar{\Delta}$ from $\theta(v)$ to $\theta(\sstar)$
to follow the arcs of $\theta \circ \gamma$ outside the interior
of $\PP_k$--islands, and to follow the geodesic path in the tree
$\theta(\partial I)$ whenever $\gamma$ enters a $\PP_k$--island
$I$.

It will be important (in \emph{Case~2} below) that the arcs of
$\bar{\gamma}$ defined as geodesics in the images of the
boundaries of $\PP_k$--islands follow edge--paths labelled by
reduced words in $\set{t, t^{-1}}^{\sstar}$. That is, they must not
traverse a pair of edges labelled by $tt^{-1}$ or
$t^{-1}t$.  The gluing $\theta$ involved choices
that, according to the following lemma (illustrated in Figure~\ref{collapse}), we can exploit to ensure
that the paths $\bar{\gamma}$ satisfy the conditions we require.

\begin{lemma}
The gluing map $\theta:\Delta \onto \bar{\Delta}$ can be chosen so
as to satisfy the following.  
Suppose $\gamma_{u,v}$ is an edge--path in $T$ from a vertex $u$ to a vertex $v$ and is an initial segment of some geodesic in $T$ from $u$  to $\sstar$.  Assume, further, that $\gamma_{u,v}$  lies in some  $\PP_k$--island $I$ and  meets $\partial I$ at  $u$ and $v$ and nowhere else.   Define $\delta_{u,v}$ to be the geodesic in the tree
$\theta(\partial I)$ from $\theta(u)$ to $\theta(v)$. Then along
$\delta_{u,v}$ we read a reduced word in $\set{t,
t^{-1}}^{\sstar}$.
\end{lemma}

\begin{figure}[ht]
\psfrag{t}{$t$}
\psfrag{u}{$\theta(u)$}
\psfrag{v}{$\theta(v)$}
\psfrag{I}{$I$}
\psfrag{s}{$\sstar$}
\psfrag{delta}{$\delta_{u,v}$}
\psfrag{gamma}{$\gamma_{u,v}$}
\psfrag{theta}{$\theta(\partial I)$}
\centerline{\epsfig{file=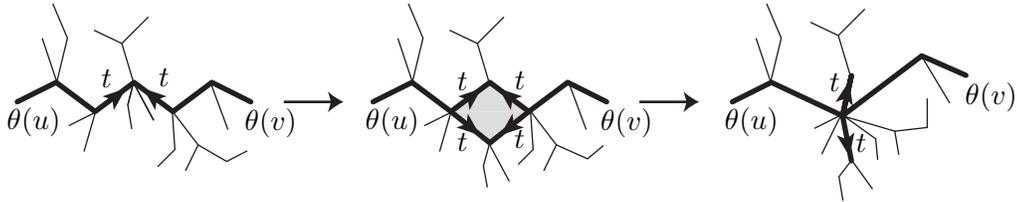}} \caption{A diamond move.}
\label{diamond}
\end{figure}

\Proof Start with any choice of gluing map $\theta$.  Suppose
that $I$ is a $\PP_k$--island in $\Delta$ and that $u, v \in
\partial I$ are as
per the lemma, but that $\delta_{u,v}$ follows a word in $\set{t,
t^{-1}}^{\sstar}$ that incudes an inverse pair, $tt^{-1}$ or
$t^{-1}t$.  Then perform a diamond move as illustrated in
Figure~\ref{diamond} to remove the inverse pair.  This, in effect,
amounts to changing the choice of $\theta$.

Suppose that $\hat{u},\hat{v} \in \partial I$ are another pair as
per the lemma. We claim that when we do the diamond move to remove
a pair of edges from $\delta_{u,v}$ the effect, if any, on the
word one reads along $\delta_{\hat{u},\hat{v}}$ is also the
removal of an inverse pair. Consider the ways in which
$\delta_{\hat{u},\hat{v}}$ could meet the pair of edges on which
the diamond move is to be performed. The danger is that an inverse
pair might be inserted into the word along
$\delta_{\hat{u},\hat{v}}$.  But this could only happen when
$\delta_{\hat{u},\hat{v}}$ \emph{crosses} $\delta_{u,v}$ at the
vertex between the two edges, and such \emph{crossing} is
impossible because $\gamma_{u,v}$ and $\gamma_{\hat{u},\hat{v}}$
are both part of geodesic arcs based at $\sstar$ in the tree $T$.

There are only finitely many pairs $u, v$ on $\partial I$. So
after a finite number of diamond moves, $\theta$ is transformed to a
map that satisfies the requirements of the lemma. \qed

\bs


Returning once more to the proof of the theorem, 
let $\bar{v}$ be a vertex in $\bar{\Delta}$ for which
$\rho(\bar{v}, \bar{\sstar}) = \IDiam(\bar{\Delta})$ and let $v$ be
a vertex of $\Delta$ such that $\theta(v) = \bar{v}$. Let $\mu$ be
the geodesic in $T$ from $v$ to $\sstar$.  And let $\bar{\mu}$ be
the edge--path from $\bar{v}$ to $\bar{\sstar}$, obtained by
connecting up the images under $\theta$ of the portions of $\mu$
not in the interior of $\PP_k$--islands $I$ with geodesic
edge--paths in the trees $\theta(I)$.  So $\bar{\mu}$ is a
concatenation of two types of geodesic arcs: those that run
through the image under $\theta$ of the boundary of some
$\PP_k$--island (call these \emph{island arcs}), and those arcs
from $\mu$.

\ms We fix $\alpha \geq \beta >0$ and examine the following two
cases.
\begin{itemize}
\item [\emph{Case 1\emph{)}}]   The island--arcs all have length at
most $n^{\alpha}$.

\ss In place of each island--arc in $\bar{\mu}$, we find an arc in
$\mu$ of length at least one (in fact, this is a crude lower
bound) because the word along each such arc in $\mu$ equals some
non--zero power of $t$ in $\Psi_{k,m}$.  So $\ell(\mu) \geq
\ell(\bar{\mu})/n^{\alpha} \geq C\, n^{m-\alpha}$ and therefore
$\IDiam(\Delta) \succeq n^{m-\alpha}$, because $\mu$ is a geodesic
in $\Delta$, based at $\sstar$.

\item [\emph{Case 2\emph{)}}] Some island--arc through
$\theta(\partial I)$, where $I$ is some $\PP_k$--island in
$\Delta$, has length more than $n^{\alpha}$.

\ss Let $\gamma$ be the corresponding subarc of $\mu$ through $I$.
If $\ell(\gamma) \geq n^{\beta}$ then it follows immediately that
$\IDiam(\Delta)  \geq \ell(\mu) \geq n^{\beta}$.

Assume that $\ell(\gamma) < n^{\beta}$. Then $\gamma$
divides $I$ into two subdiagrams each of which has boundary
circuit made up of $\gamma$ together with a portion of $\partial
\Delta$ around which we read a word in $\set{t, t^{-1}}^{\sstar}$
with exponent sum more than $n^\alpha$.  Applying
Proposition~\ref{shortcuts are fat} to either of these two
subdiagrams we learn that the intrinsic diameter of $I$ is at least a constant times
$n^{\alpha}/(1+n^{\beta})$.
\end{itemize}

\ni So $$\IDiam(\Delta) \ \succeq \ \min \set{ n^{m-\alpha},
n^{\beta}, \frac{n^{\alpha}}{1+n^{\beta}}},$$ and taking $\alpha=
2m/3$ and $\beta = m/3$ we finally have our result. \qed

\section{Amalgams and retractions} \label{Amalgams and retractions}

The main results in this section concern amalgams (that is, free
products with amalgamation). But first we give some technical results on  cutting, gluing
and collapsing operations on 2--complexes. We perform
such operations on van~Kampen diagrams and our concern is
that planarity should not be lost.

The proof of the following lemma is straight--forward and we omit it.

\begin{figure}[ht]
\psfrag{e_1}{$e_1$}
\psfrag{e_2}{$e_2$}
\psfrag{rho}{$\rho$}
\centerline{\epsfig{file=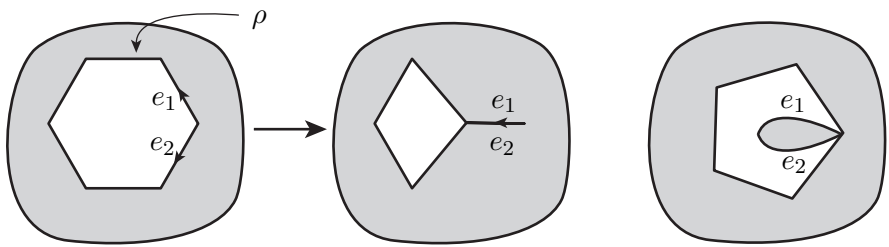}} \caption{Identification of adjacent edges.} \label{identifying}
\end{figure}

\begin{lemma} \label{collapse cut glue}
Suppose that $\Delta$ is a finite, combinatorial complex embedded
in the plane $\mathbb{E}^2$
and that $\rho$ is the \emph{(}not
necessarily embedded\emph{)} edge--circuit in $\partial \Delta$
around the boundary of the closure $\bar{C}$ of a component $C$ of
$\mathbb{E}^2-\Delta$ for which $\bar{C}$ is compact.
\begin{enumerate}
\item If $\Upsilon$ is a topological 2--disc combinatorial
complex with $\ell(\partial \Upsilon) = \ell(\rho)$, then gluing
$\Upsilon$ to $\Delta$ by identifying the boundary circuit of
$\Upsilon$ with $\rho$ produces a planar
2--complex. %
\item If $\rho$ is simple and $\Upsilon$ is a planar contractible
combinatorial 2--complex with $\ell(\partial \Upsilon) =
\ell(\rho)$, then gluing $\Upsilon$ to $\Delta$ by identifying the
boundary circuit of $\Upsilon$ with $\rho$ gives a planar
2--complex. %
\item Identifying two adjacent edges $e_1$ and $e_2$ in $\rho$, as
illustrated in Figure~\ref{identifying}, produces a planar
2--complex unless $e_1$ and $e_2$ together comprise the boundary of
a subdiagram of $\Delta$ \emph{(}as in the rightmost diagram of
the figure\emph{)}.
\end{enumerate}
\end{lemma}

A \emph{singular combinatorial map} $\Theta$ from one complex to
another is a continuous map in which every closed $n$--cell $e^n$
is either mapped homeomorphically onto an $n$--cell or is mapped
onto $\Theta(\partial e^n)$.

We leave the proof of the following technical lemma to the reader.

\begin{lemma} \label{distance decreasing}
Let $\QQ= \langle \AA \mid \RR \rangle$ and $\QQ_0= \langle \AA_0
\mid \RR_0 \rangle$ be finite presentations.  Suppose $\theta$ is
a map $\AA \to {\AA_0}^{\pm 1} \cup \set{1}$ and let $\bar{\theta}
: (\AA^{\pm 1})^\sstar \to ({\AA_0}^{\pm 1})^\sstar$ be the
extension of $\theta$ defined by
$$\bar{\theta}({a_1}^{\varepsilon_1}{a_2}^{\varepsilon_2}\ldots{a_r}^{\varepsilon_r})
={\theta(a_1)}^{\varepsilon_1}{\theta(a_2)}^{\varepsilon_2}\ldots{\theta(a_r)}^{\varepsilon_r}.$$
Suppose that for all $r \in \RR$, the word $\bar{\theta}(r)$ is
either freely reducible or has a cyclic conjugate in ${\RR_0}^{\pm
1}$.

If $\Delta$ is a $\QQ$--van~Kampen diagram for a word $w$, then
there is a singular combinatorial map $\Theta : \Delta \to
\bar{\Delta}$ to a $\QQ_0$--van~Kampen diagram for
$\bar{\theta}(w)$ that is distance decreasing with respect to the
path metrics on $\Delta^{(1)}$ and $\bar{\Delta}^{(1)}$. Moreover,
suppose $e$ is an edge of $\Delta$ labelled by $a$ and $\Theta(e)$
is not a single vertex; if $\theta(a) \in \AA$ then $\Theta$
preserves the orientation of $e$ and $\Theta(e)$ is labelled by
$\theta(a)$, and if $\theta(a) \in \AA^{-1}$ then $\Theta$
reverses the orientation of $e$ and $\Theta(e)$ is labelled by
$\theta(a)^{-1}$. \end{lemma}

\begin{cor} \label{distance decreasing cor}
In addition to the hypotheses of Lemma~\ref{distance decreasing}
assume that $\QQ_0$ is a subpresentation of $\QQ$.  If $w$ is a
null--homotopic word in $\QQ_0$, and if $\bar{\theta}(w)=w$, then
there is a $\QQ_0$--van~Kampen diagram for $w$ that is of minimal
intrinsic diameter \emph{(}or radius\emph{)} amongst all
$\QQ$--van~Kampen diagrams for $w$.
\end{cor}


In the following proposition, the hypothesis that no cyclic conjugate of a word in ${\RR_1}^{\pm1}$ is in
${\RR_2}^{\pm1}$ is mild. It could only fail for freely reducible words
in $\set{t^{\pm 1}}\*$. It allows the 2--cells of a $\PP$--van~Kampen
$\Delta$ to be partitioned into $\RR_1$--cells and
$\RR_2$--cells; that is, 2--cells that have boundary words in
$\RR_1$ or in $\RR_2$, respectively. Note that whenever an
$\RR_1$--cell shares an edge with an $\RR_2$--cell, that
edge is labelled by $t$.

\begin{prop} \label{islands}
Let $\langle \AA_1 \mid \RR_1 \rangle$ and $\langle \AA_2 \mid
\RR_2 \rangle$ be presentations of groups $A_1$ and $A_2$, such that 
$\AA_1 \cap \AA_2 = \set{t}$, where $t$ has infinite order
in both $A_1$ and $A_2$.  The amalgam $A_1 \ast_{\langle t
\rangle} A_2$ has presentation $\PP:=\langle \AA_1 \cup \AA_2 \mid
\RR_1 \cup \RR_2 \rangle$.  Assume that no cyclic conjugate of a word in ${\RR_1}^{\pm1}$ is in
${\RR_2}^{\pm1}$.

Suppose that $A_2$ retracts to
$\langle t \rangle$ via a homomorphism $\phi$ that maps $t$ to $t$
and maps all other $a \in \AA_2$ to $1$ or $t^{\pm 1}$.
\begin{enumerate}
\item Suppose $w_1$ is a null--homotopic word in
$\left({\AA_1}^{\pm1}\right)^{\sstar}$. Then $w_1$ has an $\langle
\AA_1 \mid \RR_1 \rangle$--van~Kampen diagram that is of minimal
intrinsic diameter \emph{(}or radius\emph{)} amongst all $\PP$--van~Kampen
diagrams for
$w_1$. %
\item  Suppose $w_2$ is a null--homotopic word in
$\left({\AA_2}^{\pm1}\right)^{\sstar}$. Then $w_2$ has a minimal
intrinsic diameter $\PP$--van~Kampen diagram $\hat{\Delta}$ such
that whenever $\gamma$ is a simple edge--circuit in $\hat{\Delta}$
around which we read a word in
$\left({\AA_1}^{\pm1}\right)^{\sstar}$, the subdiagram it bounds is
an $\langle \AA_1 \mid \RR_1 \rangle$--van~Kampen diagram.
\end{enumerate}
\end{prop}

\Proof The first part is a consequence of Corollary~\ref{distance
decreasing cor}. For the second part we take a $\PP$--van~Kampen
diagram $\Delta$ for $w_2$ of minimal intrinsic diameter and
obtain a $\PP$--van~Kampen diagram $\hat{\Delta}$ for $w_2$ with
the required properties by repeating the following procedure.

Suppose $\rho$ is a simple edge--circuit in $\Delta$ around which
we read a word in $\left({\AA_1}^{\pm1}\right)^{\sstar}$.  Let
$w_{\rho}$ be the word we read around $\rho$, and let
$\Delta_{\rho}$ be the $\PP$--van~Kampen subdiagram bounded by
$\rho$.  Glue the $\langle \AA_1 \mid \RR_1 \rangle$--van~Kampen
diagram $\bar{\Delta_{\rho}}$ for $w_{\rho}$ supplied by
Lemma~\ref{distance decreasing}, in place of $\Delta_{\rho}$. This
does not increase intrinsic diameter because the map $\Theta:
\Delta_{\rho} \to \bar{\Delta_{\rho}}$ of Lemma~\ref{distance
decreasing} is distance decreasing. Also the gluing cannot destroy
planarity because $\rho$ is simple --- see Lemma~\ref{collapse cut
glue} (2). \qed

\section{Quasi--isometry invariance}  \label{qi invariance}

A special case of the following theorem is that for any two finite presentations $\PP_1, \PP_2$ of the same group, $\IDiam_{\PP_1} \simeq \IDiam_{\PP_2}$ and $\EDiam_{\PP_1} \simeq \EDiam_{\PP_2}$.  

\begin{thm} \label{qi invariance thm} If $\PP_1$ and $\PP_2$ be finite
presentations for quasi--isometric groups $\G_1$ and $\G_2$, respectively, 
then $\IDiam_{\PP_1} \simeq \IDiam_{\PP_2}$ and 
$\EDiam_{\PP_1} \simeq \EDiam_{\PP_2}$. 
\end{thm}

\Proof  The theorem is proved by keeping track of diameters as one
follows the standard proof that finite presentability is a 
quasi--isometry invariant \cite[page 143]{BrH}. The first quantified version
 of this proof (in the
context of Dehn functions) appeared in \cite{Alonso}.

Fix word metrics $d_1$ and $d_2$  for $\G_1$ and $\G_2$, respectively,
and quasi--isometries $f: (\Gamma_1,d_1) \to (\Gamma_2,d_2)$ 
and $g: (\Gamma_2,d_2) \to (\Gamma_1,d_1)$ with
constants $\lambda \geq 1$ and $\mu \geq 0$, such that for all $u_1,v_1 \in \Gamma_1$ and all $u_2,v_2 \in \Gamma_2$, 
\begin{eqnarray}
  \frac{1}{\lambda} d_1(u_1,v_1) - \mu  \  \leq \   d_2(f(u_1),f(v_1)) & \leq &  \lambda d_1(u_1,v_1) + \mu, \label{f} \\ 
    \frac{1}{\lambda} d_2(u_2,v_2) - \mu  \  \leq \   d_1(g(u_2),g(v_2)) & \leq &  \lambda d_2(u_2,v_2) + \mu,  \label{g} \\ 
d(u_1, g \circ f (u_1) ),  \   d(u_2, f \circ g (u_2) )  & \leq & \mu. \label{comp}
\end{eqnarray}

Suppose $\rho_2$ is an edge--circuit in the Cayley graph
$\Cay^1(\Gamma_2, \AA_2)$, visiting vertices $v_0, v_1, \dots, v_n=v_0$
in order. Consider the circuit $\rho_1$ in $\Cay^1(\Gamma_1, \AA_1)$ 
obtained by joining  $g(v_0), g(v_1), \dots, g(v_n)$ by geodesics of length at most $\lambda+ \mu$, using (\ref{g}).  Fill $\rho_1$ with a minimal--intrinsic--diameter van~Kampen diagram $\pi_1: \Delta_1 \to \Cay^2(\PP_1)$.  Extend $f \circ \left( \pi_1\restricted{{\Delta_1}^{(0)}} \right)$ by joining the images of adjacent vertices by geodesics in $\Cay^1(\Gamma_2, \AA_2)$, each of length at  most $\lambda+ \mu$ by (\ref{f}), to give a combinatorial map $\pi: \Delta^{(1)} \to \Cay^1(\Gamma_2, \AA_2)$ from the 1--skeleton of a diagram $\Delta$ obtained  by subdividing each of the edges of $\Delta_1$.   Each 2--cell in $\Delta$ has boundary length at most $L_1:=(\lambda+ \mu)\max_{r \in \RR_1} \ell(r)$.  Extend $\pi$ to a map $\pi'_2: {\Delta_2}^{(1)} \to \Cay^1(\Gamma_2, \AA_2)$ filling $\rho_2$ by joining each vertex $u_2$ on $\rho_2$ to $f \circ g (u_2)$ on $\pi(\partial \Delta)$ by a geodesic, which  has length at most $\mu$ by (\ref{comp}).    So $\Delta_2$ is obtained from $\Delta$ by attaching a collar of $\ell(\rho_2)$ 2--cells around its boundary.  Adjacent vertices in  $\Cay^1(\Gamma_2, \AA_2)$ are mapped by $g$ to vertices at most $\lambda + \mu$ apart by (\ref{g}) and then by $f$ to vertices at most $\lambda(\lambda + \mu) + \mu$ by (\ref{f}).  So the lengths of the boundaries of the 2--cells in the collar are each at most $L_2:=\lambda^2 + \lambda \mu + 3 \mu +1$.  

It follows that if we define $\RR_2$ to be the set of all null--homotopic words in $({\AA_2}^{\pm 1})^{\ast}$ of length at most $L:=\max\set{L_1, L_2}$ then $\pi_2$ extends to a van~Kampen diagram $\pi_2:\Delta_2 \to \Cay^2(\PP_2)$ filling $\rho_2$, where $\PP_2 = \langle \AA_2 \mid \RR_2 \rangle$.  So $\PP_2$ is a finite presentation for $\Gamma_2$.  Now, $\ell(\partial \Delta_1) \leq (\lambda + \mu) \ell(\rho_2)$ and so $\IDiam(\Delta_1)$ is at most $\IDiam_{\PP_1}( (\lambda + \mu) \ell(\rho_2))$.  By (\ref{f}) we can multiply this by $(\lambda + \mu)$ to get an upper bound on the intrinsic diameter of $\Delta$.  Adding a further $2\mu$ for the collar, we get 
\begin{equation}
\IDiam(\Delta_2) \ \leq  \ (\lambda + \mu) \, \IDiam_{\PP_1}( (\lambda + \mu) \ell(\rho_2)) + 2\mu, \label{IDiam estimate}
\end{equation}
which establishes $\IDiam_{\PP_2} \preceq \IDiam_{\PP_1}$ for this particular $\PP_2$.  

However, the theorem concerns \emph{arbitrary} $\RR_2$ for which $\PP_2 = \langle \AA_2 \mid \RR_2 \rangle$ is a finite presentation for $\Gamma_2$.  The boundary of each 2--cell of $\Delta_2$ is mapped by $\pi_2$ to an edge--circuit in $\Cay^1(\Gamma_2, \AA_2)$ of length at most $L$.  So, to make $\Delta_2$ into a van~Kampen $\hat{\Delta}_2$ diagram over $\PP_2$, we fill each of its 2--cells with a (possibly singular) van~Kampen diagram over $\PP_2$.  But, a technical concern here is that gluing a \emph{singular} 2--disc diagram along a non--embedded edge--circuit of $\Delta_2$ could destroy planarity.  The way we deal with this is to fill the 2--cells of $\Delta_2$ one at a time.  And, on each occasion, if the 2--cell $C$ to be filled has non--embedded boundary circuit then we discard all the edges inside the simple edge--circuit $\sigma$ in $\partial C$ such that no edge of $\partial C$ is outside $\sigma$, and then we fill $\sigma$.   

Discarding the edges inside all such $\sigma$ does not stop the estimate (\ref{IDiam estimate}) holding.  Adding $\IDiam_{\PP_2}(L)$ to account for each of the fillings gives an upper bound on $\IDiam_{\PP_2}(\hat{\Delta}_2)$ and so $\IDiam_{\PP_2} \preceq \IDiam_{\PP_1}$.  Interchanging the roles of $\PP_1$ and $\PP_2$, we immediately deduce that $\IDiam_{\PP_1} \preceq \IDiam_{\PP_2}$ and so we have $\IDiam_{\PP_1} \simeq \IDiam_{\PP_2}$, as required.  

That $\EDiam_{\PP_1} \simeq \EDiam_{\PP_2}$ can be proved the same way, except we take $\Delta_1$ to be a minimal--extrinsic--diameter filling of $\rho_1$, and then by (\ref{f}) 
$$\EDiam(\Delta_2) \ \leq  \ \lambda \, \EDiam_{\PP_1}( (\lambda + \mu) \ell(\rho_2)) + \mu + 2L,$$
and adding a further constant $\EDiam_{\PP_2}(L)$ gives an upper bound on $\EDiam_{\PP_2}(\hat{\Delta}_2)$.  
\qed

\bibliographystyle{plain}
\bibliography{bibli}

\small{ \ni \textsc{Martin R.\ Bridson} \rule{0mm}{6mm} \\
Mathematical Institute, 24--29 St Giles', Oxford, OX1 3LB,
UK \\ \texttt{bridson@maths.ac.uk, \
http:/\!/people.maths.ox.ac.uk/$\sim$bridson/}

\ni  \textsc{Timothy R.\ Riley} \rule{0mm}{6mm} \\
Department of Mathematics, University Walk, Bristol, BS8 1TW, UK \\ \texttt{tim.riley@bristol.ac.uk, \
http:/\!/www.maths.bris.ac.uk/$\sim$matrr/ } }

\end{document}